\documentclass[12pt,thmsa]{article}
\usepackage{amssymb}
\usepackage{sw20lart}

\input tcilatex
\begin{document}

\author{Francis OGER}
\title{\textbf{Equivalence \'{e}l\'{e}mentaire entre pavages}}
\date{}

\begin{center}
\textbf{TILINGS AND ASSOCIATED RELATIONAL STRUCTURES}

Francis OGER
\end{center}

\bigskip

\bigskip

\bigskip

\noindent \textbf{Abstract.} In the present paper, as we did previously in
[7], we investigate the relations between the geometric properties of
tilings and the algebraic properties of associated relational structures.
Our study is motivated by the existence of aperiodic tiling systems. In [7],
we considered tilings of the euclidean spaces $%
\mathbb{R}
^{k}$, and isomorphism was defined up to translation. Here, we consider,
more generally, tilings of a metric space, and isomorphism is defined modulo
an arbitrary group of isometries.

In Section 1, we define the relational structures associated to tilings. The
results of Section 2 concern local isomorphism, the extraction preorder and
the characterization of relational structures which can be represented by
tilings of some given type.

In Section 3, we show that the notions of periodicity and invariance through
a translation, defined for tilings of the euclidean spaces $%
\mathbb{R}
^{k}$, can be generalized, with appropriate hypotheses, to relational
structures, and in particular to tilings of noneuclidean spaces.

The results of Section 4 are valid for uniformly locally finite relational
structures, and in particular tilings, which satisfy the local isomorphism
property. We characterize among such structures those which are locally
isomorphic to structures without nontrivial automorphism. We show that, in
an euclidean space $%
\mathbb{R}
^{k}$, this property is true for a tiling which satisfies the local
isomorphism property if and only if it is not invariant through a nontrivial
translation. We illustrate these results with examples, some of them
concerning aperiodic tiling systems of euclidean or noneuclidean
spaces.\bigskip

\noindent MSC: 52C23 (05B45, 52C22)

\noindent Keywords: Aperiodic tiling systems; Locally isomorphic tilings;
Local isomorphism property; Metric Space; Isometry; Noneuclidean; Hyperbolic
plane; Periodic; Translation; Rigid.\bigskip

\bigskip

In Section 1, we consider finite systems\ consisting of prototiles in a
metric space $E$ and local rules for assembling tiles which are equivalent
to the prototiles modulo a group $G$ of isometries of $E$. With appropriate
hypotheses, we associate to each such system $\Delta $ a finite relational
language $\mathcal{L}_{\Delta }$\ such that any $\Delta $-tiling can be
interpreted as a $\mathcal{L}_{\Delta }$-structure.

In [7], we did the same thing in a particular case: $E$ was an euclidean
space $%
\mathbb{R}
^{k}$\ and $G$ consisted of the translations of $E$. In that case, $\mathcal{%
L}_{\Delta }$ can be defined by considering all the possible configurations
of two adjacent tiles. In the present situation, we have to consider the
possible configurations for the patch consisting of all tiles within some
given \textquotedblleft distance\textquotedblright\ of one of them; that
\textquotedblleft distance\textquotedblright\ is depending on $E$, $G$ and $%
\Delta $.

In the following sections, we prove various properties of $\Delta $-tilings,
some of them classical and others new, concerning local isomorphism,
invariance through a translation... Classically, such results are shown by
considering the geometrical and topological properties of the space and the
tilings. The proofs given here are obtained by considering the algebraic
properties of the associated $\mathcal{L}_{\Delta }$-structures. In that
way, we prove the results for larger classes of tilings, and we show similar
properties for relational structures which are not represented by tilings
(see Section 4, Example 4).

It is also natural to wonder, for a given system $\Delta $, which $\mathcal{L%
}_{\Delta }$-structures can be represented by $\Delta $-tilings.
Characterizations of such structures are given in Section 2.\bigskip

\textbf{1. Tilings and associated relational structures}.\bigskip

In [7], we considered relational structures associated to tilings of the
euclidean spaces $%
\mathbb{R}
^{k}$. The isomorphism of tiles, patches, tilings... was defined up to
translation.

Actually, much more various situations have been investigated. For instance,
the isomorphism of tiles, patches, tilings... can be defined modulo a group
of isometries of $%
\mathbb{R}
^{k}$ which contains symmetries and/or rotations. Also, tilings of
noneuclidean spaces and tilings with partially overlapping tiles are
considered.

The following facts appear to be common to all cases:

\noindent - any tiling is a covering of a metric space $(E,\delta )$ by
tiles which are obtained from a finite set of closed bounded prototiles by
applying isometries belonging to some specified group;

\noindent - up to isometries in the group, only finitely many configurations
of some given size can appear in a tiling;

\noindent - the work of a tiler can be described as follows: first he puts a
small cluster of tiles somewhere in the space, then he gradually increases
the patch by adding new tiles; at each step, he applies the same finite set
of local rules, which only leave finitely many possibilities for adding new
tiles.

In order to formalize these facts, we introduce some definitions and
notations. We consider a metric space $(E,\delta )$ and a group $G$ of
bijective isometries of $E$. We do not suppose that $(E,\delta )$ is
homogeneous, or that local isometries of $E$ can be extended to global ones.

We call a \emph{tile} any closed bounded subset $T$ of $E$. We do not
suppose $T$ connected or $T$\ equal to the closure of its interior. For any
tiles $T,T^{\prime }$ (resp. any sets of tiles $\mathcal{E},\mathcal{E}%
^{\prime }$), we say that $T$ and $T^{\prime }$ (resp. $\mathcal{E}$ and $%
\mathcal{E}^{\prime }$) are \emph{isomorphic} if there exists $\sigma \in G$%
\ such that $T\sigma =T^{\prime }$ (resp. $\mathcal{E}\sigma =\mathcal{E}%
^{\prime }$). We define isomorphism in the same way for the pairs $(S,T)$
and the pairs $(\mathcal{E},T)$ with $S,T$ tiles and $\mathcal{E}$\ a set of
tiles.\bigskip

\noindent \textbf{Remark.} Sometimes, tiles with drawings are also
considered. A \emph{drawing} is a map from a tile $T$ to a finite set $%
\Omega $, whose elements are called \emph{colours}.\ In that case, the
homomorphisms that we consider must respect colours; moreover, in the
definitions of configuration, tiling and patch given below, any point in the
intersection of two tiles must have the same colour in each of them. The
results of the present paper are proved in the same way for tiles with
drawings.\bigskip

For each set $\mathcal{E}$\ of tiles and each $T\in \mathcal{E}$, we\ define
inductively the subsets $\mathcal{B}_{n}^{\mathcal{E}}(T)$\ with $\mathcal{B}%
_{0}^{\mathcal{E}}(T)=\{T\}$ and, for each $n\in 
\mathbb{N}
$,

\noindent $\mathcal{B}_{n+1}^{\mathcal{E}}(T)=\{U\in \mathcal{E\mid }$ there
exists $V\in \mathcal{B}_{n}^{\mathcal{E}}(T)$\ such that $U\cap V\neq
\emptyset \}$.

\noindent For any tiles $T,T^{\prime }$ and any sets of tiles $\mathcal{E},%
\mathcal{E}^{\prime }$, we write $(\mathcal{E},T)\leq (\mathcal{E}^{\prime
},T^{\prime })$ if there exists $\sigma \in G$\ such that $\mathcal{E}\sigma
\subset \mathcal{E}^{\prime }$ and $T\sigma =T^{\prime }$.

We call a \emph{tiling} any covering $\mathcal{E}$\ of $E$ by possibly
overlapping tiles such that, for each $r\in 
\mathbb{N}
^{\ast }$:

\noindent 1) for each $T\in \mathcal{E}$, $\mathcal{B}_{r}^{\mathcal{E}}(T)$
is finite and $\mathcal{B}_{r-1}^{\mathcal{E}}(T)$ is contained in the
interior of the union of the tiles of $\mathcal{B}_{r}^{\mathcal{E}}(T)$;

\noindent 2) for any $S,T\in \mathcal{E}$, if $(\mathcal{B}_{r}^{\mathcal{E}%
}(S),S)\leq (\mathcal{B}_{r}^{\mathcal{E}}(T),T)$, then $(\mathcal{B}_{r}^{%
\mathcal{E}}(S),S)\cong (\mathcal{B}_{r}^{\mathcal{E}}(T),T)$;

\noindent 3) the pairs $(\mathcal{B}_{r}^{\mathcal{E}}(T),T)$ for $T\in 
\mathcal{E}$ fall in finitely many isomorphism classes.

If 1) is true for $r=1$, then it is true for each $r\in 
\mathbb{N}
^{\ast }$. If 1) is true and if 2) is true for $r=1$, then, for any $S,T\in 
\mathcal{E}$, each $\sigma \in G$\ such that $\mathcal{B}_{1}^{\mathcal{E}%
}(S)\sigma \subset \mathcal{B}_{1}^{\mathcal{E}}(T)$ and $S\sigma =T$\ is an
isomorphism from $(\mathcal{B}_{1}^{\mathcal{E}}(S),S)$ to $(\mathcal{B}%
_{1}^{\mathcal{E}}(T),T)$. By induction, it follows that, for each $r\in 
\mathbb{N}
^{\ast }$ and any $S,T\in \mathcal{E}$, each $\sigma \in G$\ such that $%
\mathcal{B}_{r}^{\mathcal{E}}(S)\sigma \subset \mathcal{B}_{r}^{\mathcal{E}%
}(T)$ and $S\sigma =T$ is an isomorphism from $(\mathcal{B}_{r}^{\mathcal{E}%
}(S),S)$ to $(\mathcal{B}_{r}^{\mathcal{E}}(T),T)$. Consequently, 2) is true
for each $r\in 
\mathbb{N}
^{\ast }$. We note that 2) is necessarily true if 1) is true and if the
tiles of $\mathcal{E}$ are equal to the closure of their interior and
nonoverlapping.

It is natural to wonder whether 3) is also true for each $r\in 
\mathbb{N}
^{\ast }$ if it is true for $r=1$. If $(E,\delta )$ is an euclidean space $%
\mathbb{R}
^{k}$ and if $G$ consists of the translations of $%
\mathbb{R}
^{k}$, then any pair $(S,T)\in \mathcal{E}\times \mathcal{E}$ is completely
determined by $S$ and the isomorphism class of $(S,T)$. Consequently, if the
pairs $(S,T)\in \mathcal{E}\times \mathcal{E}$ with $S\cap T\neq \emptyset $
fall in finitely many isomorphism classes, then, for each $n\in 
\mathbb{N}
^{\ast }$, the pairs $(\mathcal{B}_{n}^{\mathcal{E}}(T),T)$ for $T\in 
\mathcal{E}$ fall in finitely many isomorphism classes. The following
example shows that the situation is different if $G$ consists of the
isometries of $%
\mathbb{R}
^{k}$:\bigskip

\noindent \textbf{Example 1.} Let $(E,\delta )$ be the euclidean space $%
\mathbb{R}
^{2}$ and let $G$ consist of the isometries of $(E,\delta )$. Consider the
coverings of $%
\mathbb{R}
^{2}$ by nonoverlapping tiles isomorphic to the following prototiles:

\noindent $T_{0}=\{(x,y)\in 
\mathbb{R}
^{2}\mid x\geq 0$ and $x^{2}+y^{2}\leq 1\}$,

\noindent $T_{k}=\{(x,y)\in 
\mathbb{R}
^{2}\mid k^{2}\leq x^{2}+y^{2}\leq (k+1)^{2}\}$ for $1\leq k\leq 2n+1$,

\noindent $T_{2n+2}=\{(x,y)\in 
\mathbb{R}
^{2}\mid x^{2}+y^{2}\geq (2n+2)^{2}$ and $\sup (\left\vert x\right\vert
,\left\vert y\right\vert )\leq 2n+3\}$,

\noindent with any two copies of $T_{2n+2}$ having one edge in common if
they have more than one common point.

Each such covering consists of squares with sides of length $4n+6$, each of
them covered by two copies of $T_{0}$ and one copy of each of the tiles $%
T_{1},...,T_{2n+2}$. Moreover, in each square, the two copies of $T_{0}$ can
be rotated together arbitrarily. Consequently, in such a covering $\mathcal{E%
}$, the pairs $(\mathcal{B}_{n+1}^{\mathcal{E}}(T),T)$ for $T\in \mathcal{E}$
do not generally fall in finitely many isomorphism classes, contrary to the
pairs $(\mathcal{B}_{n}^{\mathcal{E}}(T),T)$.\bigskip

Now, we can ask the following question: Suppose that we specify a finite set
of isomorphism classes for the pairs $(\mathcal{B}_{1}^{\mathcal{E}}(T),T)$
which can appear in a covering $\mathcal{E}$\ of $E$ which satisfies 1) and
2). Is there an integer $m$ such that, if the pairs $(\mathcal{B}_{m}^{%
\mathcal{E}}(T),T)$ for $T\in \mathcal{E}$ fall in finitely many isomorphism
classes, then, for each $n\in 
\mathbb{N}
^{\ast }$, the pairs $(\mathcal{B}_{n}^{\mathcal{E}}(T),T)$ for $T\in 
\mathcal{E}$ fall in finitely many isomorphism classes.

Example 2 below shows that this property does not necessarily hold. On the
other hand, we are going to prove (see Proposition 1.10) that some natural
conditions on $(E,\delta )$ and $G$, satisfied by the euclidean spaces $%
\mathbb{R}
^{k}$, and on the pairs $(\mathcal{B}_{1}^{\mathcal{E}}(T),T)$ for $T\in 
\mathcal{E}$, are enough to make it true.\bigskip

\noindent \textbf{Example 2.} Let $E=\{(x,y,z)\in 
\mathbb{R}
^{3}\mid x^{2}+y^{2}=1\}$ be the surface of a cylinder of infinite length.
Define the distance $\delta $ by considering $E$ as a quotient of the
euclidean space $%
\mathbb{R}
^{2}$. Let $G$ consist of the isometries of $E$. Write $T_{1}=\{(x,y,z)\in
E\mid x\geq 0$ and $0\leq z\leq 1\}$ and $T_{2}=\{(x,y,z)\in E\mid 0\leq
z\leq 1\}$. For each $n\in 
\mathbb{N}
^{\ast }$, consider the coverings $\mathcal{E}$ of $E$ by nonoverlapping
tiles obtained as follows: for each $a\in 
\mathbb{Z}
$, $\{(x,y,z)\in E\mid a\leq z\leq a+1\}$ is covered by two tiles isomorphic
to $T_{1}$ if $a\in (2n+2)%
\mathbb{Z}
$\ and by one tile isomorphic to $T_{2}$ otherwise. For each $a\in (2n+2)%
\mathbb{Z}
$, the two tiles covering $\{(x,y,z)\in E\mid a\leq z\leq a+1\}$ can be
rotated together arbitrarily. Consequently, in such a covering $\mathcal{E}$%
, the pairs $(\mathcal{B}_{n+1}^{\mathcal{E}}(T),T)$ for $T\in \mathcal{E}$
do not generally fall in finitely many isomorphism classes, contrary to the
pairs $(\mathcal{B}_{n}^{\mathcal{E}}(T),T)$. Moreover, the pairs $(\mathcal{%
B}_{1}^{\mathcal{E}}(T),T)$ for all possible choices of $\mathcal{E}$ and $T$
fall in three isomorphism classes.\bigskip

Now we state the definitions, the notations and the hypotheses which are
necessary for our results. For each $x\in E$ and each $\eta \in 
\mathbb{R}
_{\geq 0}$, we write $\beta (x,\eta )=\{y\in E\mid \delta (x,y)\leq \eta \}$%
. For each nonempty bounded $S\subset E$, we consider the $\emph{radius}$ $%
\mathrm{Rad}(S)=\inf_{x\in E}(\sup_{y\in S}\delta (x,y))$ and the \emph{%
diameter} $\mathrm{Diam}(S)=\sup_{x,y\in S}\delta (x,y)$; we have $\mathrm{%
Diam}(S)\leq 2.\mathrm{Rad}(S)$.

From now on, we suppose that $E$ is connected and that $\beta (x,\eta )$ is
compact for each $x\in E$ and each $\eta \in 
\mathbb{R}
_{>0}$. Then each bounded closed subset (i.e. tile) of $E$ is compact.

For each $x\in E$ and each $\xi \in 
\mathbb{R}
_{>0}$, we define inductively $\beta (x,\xi ,n)$ with $\beta (x,\xi
,0)=\left\{ x\right\} $ and, for each $n\in 
\mathbb{N}
$, $\beta (x,\xi ,n+1)=\cup _{y\in \beta (x,\xi ,n)}\beta (y,\xi )$. We have 
$\cup _{n\in 
\mathbb{N}
}\beta (x,\xi ,n)=E$\ since $E$ is connected and $\cup _{n\in 
\mathbb{N}
}\beta (x,\xi ,n)$\ is both open and closed. For each $n\in 
\mathbb{N}
$, we write $\omega (x,\xi ,n)=\sup \{\eta \in 
\mathbb{R}
_{\geq 0}\mid \beta (x,\eta )\subset \beta (x,\xi ,n)\}$.

Now we show that $\omega (x,\xi ,n)$\ tends to infinity with $n$. Otherwise,
there exists $\eta \in 
\mathbb{R}
_{>0}$\ such that $\beta (x,\eta )\nsubseteqq \beta (x,\xi ,n)$\ for each $%
n\in 
\mathbb{N}
$. It follows that $\beta (x,\eta )$ contains a sequence $(y_{n})_{n\in 
\mathbb{N}
}$\ with $y_{n}\notin \beta (x,\xi ,n)$\ for each $n\in 
\mathbb{N}
$. As $\beta (x,\eta )$ is compact, there exists a subsequence $%
(z_{n})_{n\in 
\mathbb{N}
}$\ which converges to a point $z\in \beta (x,\eta )$. For $n$ large enough,
we have\ $\delta (z,z_{n})<\xi $, which implies that $z\notin \beta (x,\xi
,n-1)$\ since $z_{n}\notin \beta (x,\xi ,n)$. It follows that $z\notin \cup
_{n\in 
\mathbb{N}
}\beta (x,\xi ,n)$, contrary to what we proved just above.

We say that $(E,\delta )$\ is:

\noindent -\ \emph{weakly homogeneous} if $\omega (\xi ,n)=\inf_{x\in
E}\omega (x,\xi ,n)$\ tends to infinity with $n$;

\noindent -\ \emph{geodesic}\ if any two points of $E$\ can be joined by at
least one geodesic;

\noindent -\emph{\ transitive} if, for any $x,y\in E$, there exists an
isometry $\sigma $\ such that $x\sigma =y$.

Any geodesic space is weakly homogeneous since $\omega (x,\xi ,n)=n\xi $ for
any $x,\xi ,n$. Any transitive space is also weakly homogeneous since $%
\omega (x,\xi ,n)$\ does not depend on $x$. In some of our results, we shall
suppose that $(E,\delta )$\ is weakly homogeneous or geodesic.

For each \emph{finite} set $\mathcal{E}$ of tiles, each $T\in \mathcal{E}$
and each $p\in 
\mathbb{N}
^{\ast }$, we say that $(\mathcal{E},T)$ is a $p$\emph{-configuration} if $%
\mathcal{B}_{p}^{\mathcal{E}}(T)=\mathcal{E}$ and if each $S\in \mathcal{B}%
_{p-1}^{\mathcal{E}}(T)$ is contained in the interior of the union of the
tiles of $\mathcal{B}_{1}^{\mathcal{E}}(S)$. Then $\mathcal{B}_{p-1}^{%
\mathcal{E}}(T)$ is contained in the interior of the union of the tiles of $%
\mathcal{E}$ and $(\mathcal{B}_{k}^{\mathcal{E}}(S),S)$ is a $k$%
-configuration for each $k\in \{1,...,p\}$ and each $S\in \mathcal{B}_{p-k}^{%
\mathcal{E}}(T)$.\bigskip

\noindent \textbf{Lemma 1.1.} For each $1$-configuration $(\mathcal{E},T)$,
there exists $\xi \in 
\mathbb{R}
_{>0}$ such that $\cup _{x\in T}\beta (x,\xi )$\ is contained in the union
of the tiles of $\mathcal{E}$.\bigskip

\noindent \textbf{Proof.} For each $n\in 
\mathbb{N}
^{\ast }$, consider $x_{n}\in T$ such that $\beta (x_{n},1/n)$\ is not
contained in the union of the tiles of $\mathcal{E}$. As $T$ is compact,
there exists a subsequence\ of $(x_{n})_{n\in 
\mathbb{N}
}$\ which converges to a point $x\in T$, and $x$ does not belong to the
interior of the union of the tiles of $\mathcal{E}$, whence a
contradiction.~~$\blacksquare $\bigskip

\noindent \textbf{Proposition 1.2.} Let $\mathcal{E}$\ be a nonempty set of
tiles such that the pairs $(\mathcal{B}_{1}^{\mathcal{E}}(T),T)$ for $T\in 
\mathcal{E}$ are $1$-configurations and fall in finitely many isomorphism
classes. Then

\noindent 1) $\mathcal{E}$\ covers $E$ and $\{S\in \mathcal{E}\mid S\cap
\beta (x,\eta )\neq \varnothing \}$ is finite for each $x\in E$ and each $%
\eta \in 
\mathbb{R}
_{>0}$;

\noindent 2) If $(E,\delta )$\ is weakly homogeneous, then, for each $\eta
\in 
\mathbb{R}
_{>0}$, there exists $r\in 
\mathbb{N}
^{\ast }$ such that $\delta (S,T)\leq \eta $ implies $T\in \mathcal{B}_{r}^{%
\mathcal{E}}(S)$ for any $S,T\in \mathcal{E}$.\bigskip

\noindent \textbf{Proof.} By Lemma 1.1, there exists $\xi \in 
\mathbb{R}
_{>0}$ such that $\mathcal{B}_{1}^{\mathcal{E}}(T)$\ covers $\beta (x,\xi )$%
\ for each $T\in \mathcal{E}$ and each $x\in T$. Consequently, the union of
the tiles of $\mathcal{E}$ is both open and closed, and $\mathcal{E}$\
covers $E$ since $E$ is connected.

Now suppose that there exist $x\in E$ and $\eta \in 
\mathbb{R}
_{>0}$\ such that $\{S\in \mathcal{E}\mid S\cap \beta (x,\eta )\neq
\varnothing \}$ is infinite. Consider a sequence $(y_{n})_{n\in 
\mathbb{N}
}\subset \beta (x,\eta )$\ and a sequence of distinct tiles $(T_{n})_{n\in 
\mathbb{N}
}\subset \mathcal{E}$ such that $y_{n}\in T_{n}$ for each $n\in 
\mathbb{N}
$. As $\beta (x,\eta )$\ is compact, there exists a subsequence of $%
(y_{n})_{n\in 
\mathbb{N}
}$ which converges to a point $y\in \beta (x,\eta )$. For each $T\in 
\mathcal{E}$\ such that $y\in T$, $\mathcal{B}_{2}^{\mathcal{E}}(T)$\
contains infinitely many $T_{n}$\ since $\mathcal{B}_{1}^{\mathcal{E}}(T)$\
covers $\beta (x,\xi )$, whence a contradiction.

The second part of Proposition 1.2 is true because, for each $n\in 
\mathbb{N}
$, each $S\in \mathcal{E}$ and each $x\in S$, $\mathcal{B}_{n}^{\mathcal{E}%
}(S)$\ covers $\beta (x,\omega (\xi ,n))$\ and $\mathcal{B}_{n+1}^{\mathcal{E%
}}(S)$\ contains $\{T\in \mathcal{E}\mid T\cap \beta (x,\omega (\xi ,n))\neq
\varnothing \}$.~~$\blacksquare $\bigskip

For each $p\in 
\mathbb{N}
^{\ast }$, we call a $p$\emph{-local rule} any set $\Gamma =\{(\mathcal{C}%
_{1},C_{1}),...,(\mathcal{C}_{m},C_{m})\}$\ of pairwise nonisomorphic $p$%
-configurations such that:

\noindent 1) for any $i,j\in \{1,...,m\}$, if\ $(\mathcal{B}_{1}^{\mathcal{C}%
_{i}}(C_{i}),C_{i})\leq (\mathcal{B}_{1}^{\mathcal{C}_{j}}(C_{j}),C_{j})$,
then $(\mathcal{B}_{1}^{\mathcal{C}_{i}}(C_{i}),C_{i})\cong (\mathcal{B}%
_{1}^{\mathcal{C}_{j}}(C_{j}),C_{j})$;

\noindent 2) for each $i\in \{1,...,m\}$ and each $S\in \mathcal{B}_{1}^{%
\mathcal{C}_{i}}(C_{i})$, there exists $j\in \{1,...,m\}$ such that $(%
\mathcal{B}_{p-1}^{\mathcal{C}_{i}}(S),S)\cong (\mathcal{B}_{p-1}^{\mathcal{C%
}_{j}}(C_{j}),C_{j})$.

\noindent If $\Gamma $ is a $p$-local rule, then, for each $q\in
\{1,...,p-1\}$, any representatives of the isomorphism classes of $(\mathcal{%
B}_{q}^{\mathcal{C}_{1}}(C_{1}),C_{1}),...,(\mathcal{B}_{q}^{\mathcal{C}%
_{m}}(C_{m}),C_{m})$ form a $q$-local rule.

We say that a set $\mathcal{E}$\ of tiles $\emph{satisfies}$ $\Gamma $\ if,
for each $T\in \mathcal{E}$, the pair $(\mathcal{B}_{p}^{\mathcal{E}}(T),T)$
is isomorphic to one of the pairs $(\mathcal{C}_{i},C_{i})$. Each tiling
satisfies a $p$-local rule for each $p\in 
\mathbb{N}
^{\ast }$.

We say that a set $\mathcal{E}$\ of tiles is a \emph{patch} if there exist
no $\mathcal{F},\mathcal{G}\subset \mathcal{E}$\ such that $\mathcal{F}\cup 
\mathcal{G}=\mathcal{E}$\ and $S\cap T=\varnothing $\ for each $S\in 
\mathcal{F}$\ and each $T\in \mathcal{G}$.\ This property is true if and
only if, for any $S,T\in \mathcal{E}$, there exists $r\in 
\mathbb{N}
$ such that $T\in \mathcal{B}_{r}^{\mathcal{E}}(S)$.

For each finite set $\mathcal{E}$\ of tiles and each connected set $A\subset
E$, if $A$ is contained in the union of the tiles of $\mathcal{E}$, and if
each tile of $\mathcal{E}$ contains a point of $A$, then $\mathcal{E}$ is a
patch. The union of the tiles of a patch is not necessarily connected if the
tiles themselves are not connected.

For each set $\mathcal{E}$ of tiles, we write $\mathrm{I}_{\Gamma }(\mathcal{%
E})=\cup _{1\leq i\leq m}\{T\in \mathcal{E}\mid (\mathcal{B}_{p}^{\mathcal{E}%
}(T),T)\cong (\mathcal{C}_{i},C_{i})\}$. We say that $\mathcal{E}$ is a $%
\Gamma $\emph{-patch} if:

\noindent 1) for each $T\in \mathcal{E}$, there exists $i\in \{1,...,m\}$
such that $(\mathcal{B}_{p}^{\mathcal{E}}(T),T)\leq (\mathcal{C}_{i},C_{i})$.

\noindent 2)\ $\mathrm{I}_{\Gamma }(\mathcal{E})$\ contains a patch $%
\mathcal{A}$ such that $\mathcal{E}=\cup _{T\in \mathcal{A}}\mathcal{B}_{p}^{%
\mathcal{E}}(T)$.

For each set $\mathcal{E}$\ of tiles which satisfies $\Gamma $,\ each\ $T\in 
\mathcal{E}$\ and each integer $k\geq p$, $\mathcal{B}_{k}^{\mathcal{E}}(T)$
is a $\Gamma $-patch since $\mathcal{B}_{k-p}^{\mathcal{E}}(T)$\ is
connected, $\mathcal{B}_{k-p}^{\mathcal{E}}(T)\subset \mathrm{I}_{\Gamma }(%
\mathcal{B}_{k}^{\mathcal{E}}(T))$\ and $\mathcal{B}_{k}^{\mathcal{E}%
}(T)=\cup _{S\in \mathcal{B}_{k-p}^{\mathcal{E}}(T)}\mathcal{B}_{p}^{%
\mathcal{E}}(S)$.

Now, for each $q\in 
\mathbb{N}
^{\ast }$ and each $q$-local rule $\Delta =\{(\mathcal{D}_{1},D_{1}),...,(%
\mathcal{D}_{n},D_{n})\}$ such that each $\mathcal{B}_{q-1}^{\mathcal{D}%
_{i}}(D_{i})$ is fixed by no $\sigma \in G-\{\mathrm{Id}\}$, we define a
finite relational language $\mathcal{L}_{\Delta }$ such that the $\Delta $%
-patches can be represented by $\mathcal{L}_{\Delta }$-structures. As a
consequence, we show that any set\ of tiles which satisfies $\Delta $\ is a
tiling. Practically, in many examples of tilings, each $1$-configuration is
fixed by no $\sigma \in G-\{\mathrm{Id}\}$, so that we can take $\Delta $
with $q=2$.

For each $i\in \{1,...,n\}$, we write $\mathcal{D}_{i}=%
\{D_{i,1},...,D_{i,p(i)}\}$ with $D_{i,1}=D_{i}$, and we introduce a $p(i)$%
-ary relational symbol $R_{i}(u_{i,1},...,u_{i,p(i)})$. We write $\mathcal{L}%
_{\Delta }=\{R_{1},...,R_{n}\}$.

For each set $\mathcal{E}$ of tiles, we define a $\mathcal{L}_{\Delta }$%
\emph{-structure} on $\mathcal{E}$ as follows: for $1\leq i\leq n$\ and $%
T_{1},...,T_{p(i)}\in \mathcal{E}$, we write $R_{i}(T_{1},...,T_{p(i)})$\ if 
$\mathcal{B}_{q}^{\mathcal{E}}(T_{1})=\{T_{1},...,T_{p(i)}\}$ and if there
exists $\sigma \in G$\ such that $D_{i,j}\sigma =T_{j}$ for $1\leq j\leq
p(i) $. For $1\leq i,j\leq n$ and $S,S_{2},...,S_{p(i)},T_{2},...,T_{p(j)}%
\in \mathcal{E}$,\ the relations $R_{i}(S,S_{2},...,S_{p(i)})$ and $%
R_{j}(S,T_{2},...,T_{p(j)})$ imply $i=j$\ and $\{S_{2},...,S_{p(i)}\}=%
\{T_{2},...,T_{p(i)}\}$.

Any set of tiles $\mathcal{E}$ satisfies $\Delta $ if and only if, for each $%
T\in \mathcal{E}$,\ there exist $1\leq i\leq n$\ and $T_{2},...,T_{p(i)}\in 
\mathcal{E}$\ such that $R_{i}(T,T_{2},...,T_{p(i)})$. For any $\Delta $%
\textbf{-}patches $\mathcal{E}\subset \mathcal{F}$, the $\mathcal{L}_{\Delta
}$-structure defined on $\mathcal{E}$ is the restriction to $\mathcal{E}$ of
the $\mathcal{L}_{\Delta }$-structure defined on $\mathcal{F}$.

The following theorem implies that any $\Delta $\textbf{-}patch, and in
particular any set of tiles which satisfies $\Delta $, is determined up to
isomorphism by the associated $\mathcal{L}_{\Delta }$-structure:\bigskip

\noindent \textbf{Theorem 1.3.} For any $\Delta $\textbf{-}patches $\mathcal{%
E},\mathcal{F}$ and each $\mathcal{L}_{\Delta }$-homomorphism $f:\mathcal{E}%
\rightarrow \mathcal{F}$, there exists a unique $\sigma \in G$ such that $%
Sf=S\sigma $ for each $S\in \mathcal{E}$.\bigskip

\noindent \textbf{Proof.} Consider a patch $\mathcal{A}\subset \mathrm{I}%
_{\Delta }(\mathcal{E})$ such that $\mathcal{E}=\cup _{T\in \mathcal{A}}%
\mathcal{B}_{q}^{\mathcal{E}}(T)$. For each $T\in \mathcal{A}$, $f$ induces
a bijection from $\mathcal{B}_{q}^{\mathcal{E}}(T)$ to $\mathcal{B}_{q}^{%
\mathcal{F}}(Tf)$, and there exists a unique $\sigma _{T}\in G$ such that $%
S\sigma _{T}=Sf$ for each $S\in \mathcal{B}_{q}^{\mathcal{E}}(T)$.

It remains to be proved that $\sigma _{S}=\sigma _{T}$\ for any $S,T\in 
\mathcal{A}$.\ As $\mathcal{A}$ is a patch, it suffices to show it for $%
S\cap T\neq \emptyset $. Then we have $S\in \mathcal{B}_{1}^{\mathcal{E}}(T)$%
, and therefore $\mathcal{B}_{q-1}^{\mathcal{E}}(S)\subset \mathcal{B}_{q}^{%
\mathcal{E}}(T)$. Consequently, for each $U\in \mathcal{B}_{q-1}^{\mathcal{E}%
}(S)$, we have $U\sigma _{S}=Uf$ since $U\in \mathcal{B}_{q}^{\mathcal{E}%
}(S) $, and $U\sigma _{T}=Uf$ since $U\in \mathcal{B}_{q}^{\mathcal{E}}(T)$.
It follows $U\sigma _{S}=U\sigma _{T}$ for each $U\in \mathcal{B}_{q-1}^{%
\mathcal{E}}(S)$, and therefore\ $\sigma _{S}=\sigma _{T}$.~~$\blacksquare $%
\bigskip

\noindent \textbf{Corollary 1.4.} For each integer $k$, there are finitely
many isomorphism classes of $\Delta $\textbf{-}patches consisting of $k$
tiles.\bigskip

\noindent \textbf{Proof.} This result follows from Theorem 1.3 since there
are finitely many isomorphism classes of $\mathcal{L}_{\Delta }$-structures
consisting of $k$ elements.~~$\blacksquare $\bigskip

\noindent \textbf{Corollary 1.5.} For each finite $\Delta $\textbf{-}patch $%
\mathcal{E}$ and each $k\in 
\mathbb{N}
^{\ast }$, finitely many distinct $\Delta $\textbf{-}patches can be obtained
from $\mathcal{E}$ by adding $k$ new tiles.\bigskip

\noindent \textbf{Proof.} Denote by $\Omega $\ the set of all $\Delta $%
\textbf{-}patches obtained from $\mathcal{E}$ by adding $k$ new tiles. By
corollary 1.4, it suffices to show that, for each $\mathcal{F}\in \Omega $,
there exist finitely many pairs $(\mathcal{G},\sigma )$ with $\mathcal{G}\in
\Omega $\ and $\sigma :\mathcal{F}\rightarrow \mathcal{G}$\ isomorphism. But
any such pair is completely determined by the tiles $\sigma ^{-1}(T)$\ for $%
T\in \mathcal{E}$, since no $\rho \in G-\{\mathrm{Id}\}$ fixes the tiles of $%
\mathcal{E}$.~~$\blacksquare $\bigskip

\noindent \textbf{Corollary 1.6.} For each $k\in 
\mathbb{N}
$, $\{(\mathcal{B}_{k}^{\mathcal{E}}(T),T)\mid $ $\mathcal{E}$ satisfies $%
\Delta $ and $T\in \mathcal{E}\}$\ is a finite union of isomorphism
classes.\bigskip

\noindent \textbf{Proof.} For $k\leq q$, this property is true because $%
\mathcal{B}_{k}^{\mathcal{E}}(T)$\ is contained in $\mathcal{B}_{q}^{%
\mathcal{E}}(T)$ for each set $\mathcal{E}$ of tiles which satisfies $\Delta 
$ and each $T\in \mathcal{E}$. For $k\geq q+1$, it follows from Corollary
1.4 since we have a finite bound for the cardinals of the $\Delta $\textbf{-}%
patches $\mathcal{B}_{k}^{\mathcal{E}}(T)$ for $\mathcal{E}$ satisfying $%
\Delta $ and $T\in \mathcal{E}$.~~$\blacksquare $\bigskip

\noindent \textbf{Corollary 1.7.} Any set of tiles which satisfies $\Delta $
is a tiling.\bigskip

\noindent \textbf{Proof.} Any such set $\mathcal{E}$ covers $E$\ by the
first part of Proposition 1.2. Moreover, for each $k\in 
\mathbb{N}
$, $\{(\mathcal{B}_{k}^{\mathcal{E}}(T),T)\mid T\in \mathcal{E}\}$\ is a
finite union of isomorphism classes by Corollary 1.6.~~$\blacksquare $%
\bigskip

\noindent \textbf{Remark.} The following variant of Example 1 above shows
that, in Theorem 1.3 and its corollaries, it is not enough to suppose that
each $\mathcal{D}_{i}$ is fixed by no $\sigma \in G-\{\mathrm{Id}\}$.
Consider the coverings $\mathcal{E}$ of $%
\mathbb{R}
^{2}$\ obtained from the following prototiles:\ 

\noindent $T_{0}=\{(x,y)\in 
\mathbb{R}
^{2}\mid x\geq 0$ and $x^{2}+y^{2}\leq 1\}$,

\noindent $T_{k}=\{(x,y)\in 
\mathbb{R}
^{2}\mid k^{2}\leq x^{2}+y^{2}\leq (k+1)^{2}\}$ for $1\leq k\leq 2q$,

\noindent $T_{2q+1}=\{(x,y)\in 
\mathbb{R}
^{2}\mid x^{2}+y^{2}\geq (2q+1)^{2}$ and $\sup (\left\vert x\right\vert
,\left\vert y\right\vert )\leq 2q+2\}$,

\noindent with some bumps on the diameter of $T_{0}$ and on the four sides
of $T_{2q+1}$ so that $T_{0}$\ and $T_{2q+1}$ are fixed by no $\sigma \in
G-\{\mathrm{Id}\}$. Then, in each such $\mathcal{E}$, each $\mathcal{B}_{q}^{%
\mathcal{E}}(T)$ is fixed by no $\sigma \in G-\{\mathrm{Id}\}$, but the
pairs $(\mathcal{B}_{q+1}^{\mathcal{E}}(T),T)$ for $T\in \mathcal{E}$ do not
generally fall in finitely many isomorphism classes.\bigskip

We say that a group $H$\ of isometries of $E$\ is \emph{discrete} if, for
each $x\in E$, there exists $\eta \in 
\mathbb{R}
_{>0}$ such that $\beta (x,\eta )\cap \{x\sigma \mid \sigma \in H\}=\{x\}$.
This property is true if and only if $\beta (y,\eta )\cap \{x\sigma \mid
\sigma \in H\}$ is finite for each $\eta \in 
\mathbb{R}
_{>0}$ and any $x,y\in E$.\bigskip

\noindent \textbf{Proposition 1.8.} For each tiling $\mathcal{T}$ which
satisfies $\Delta $, the subgroup $H=\{\sigma \in G\mid \mathcal{T}\sigma =%
\mathcal{T}\}$\ is discrete.\bigskip

\noindent \textbf{Proof.} We show that $\beta (y,\eta )\cap \{x\sigma \mid
\sigma \in H\}$ is finite for each $\eta \in 
\mathbb{R}
_{>0}$ and any $x,y\in E$. We consider $T\in \mathcal{T}$\ such that $x\in T$%
. Any $\sigma \in H$\ such that $x\sigma \in \beta (y,\eta )$\ sends $T$ to
a tile $U$ such that $U\cap \beta (y,\eta )\neq \varnothing $ and induces a
bijection from $\mathcal{B}_{q}^{\mathcal{T}}(T)$\ to $\mathcal{B}_{q}^{%
\mathcal{T}}(U)$\ which completely determines $\sigma $. By the first part
of Proposition 1.2, $\mathcal{T}$ contains finitely many tiles $U$ such that 
$U\cap \beta (y,\eta )\neq \varnothing $. Moreover, for each such $U$, there
exist finitely many bijections from $\mathcal{B}_{q}^{\mathcal{T}}(T)$\ to $%
\mathcal{B}_{q}^{\mathcal{T}}(U)$. Consequently, $\{\sigma \in H\mid x\sigma
\in \beta (y,\eta )\}$\ is finite.~~$\blacksquare $\bigskip

Now we fix $\lambda \in 
\mathbb{R}
_{>0}$ and we state some supplementary conditions on $(E,\delta )$ and $G$
which imply that each tiling with tiles of radius $<\lambda $ satisfies a
local rule with the properties stated above. We suppose that $(E,\delta )$
is geodesic and we consider the following properties:

\noindent (CVX$\lambda $) For each $x\in E$ and each $\eta \in \left]
0,\lambda \right[ $, each geodesic which joins two points of $\beta (x,\eta
) $\ is contained in $\cup _{\zeta <\eta }\beta (x,\zeta )$;

\noindent (FIX$\lambda $)\ There exists $\mu \in 
\mathbb{R}
_{>0}$ such that:

\noindent (FIX$\lambda \mu $) for each $x\in E$, no $\sigma \in G-\{\mathrm{%
Id}\}$ fixes the points of a set $A\subset E$\ with $\beta (x,\mu )\subset
\cup _{y\in A}\beta (y,\lambda )$.

In euclidean or hyperbolic spaces of finite dimension, (CVX$\nu $)\ and (FIX$%
\nu $)\ are true for each $\nu \in 
\mathbb{R}
_{>0}$. On the other hand, for the surfaces of a sphere or a cylinder of
infinite length, they are only true for $\nu $ small enough. We do not want
to suppose from the beginning that they are true for each $\nu \in 
\mathbb{R}
_{>0}$ since, for instance, any tiling of $%
\mathbb{R}
^{2}$\ which is invariant through a nontrivial translation induces a tiling
of the surface of a cylinder of infinite length.\bigskip

\noindent \textbf{Lemma 1.9.} Suppose that $(E,\delta )$ is geodesic and
satisfies (CVX$\lambda $). Then, for each nonempty $S\subset E$ such that $%
\mathrm{Rad}(S)<\lambda $, there exists a unique $x\in E$ such that $\delta
(x,y)\leq \mathrm{Rad}(S)$\ for each $y\in S$, and we have $x\sigma =x$\ for
each isometry $\sigma $\ such that $S\sigma =S$.\bigskip

\noindent \textbf{Proof.} We only prove the first statement, since the
second one is an immediate consequence. For each $w\in E$, we have $%
\sup_{y\in S}\delta (w,y)=\sup_{y\in T}\delta (w,y)$ where $T$ is the
closure of $S$ in $E$. Consequently, we can suppose $S$ closed, and
therefore $S$ compact. The subset $A=\{w\in E\mid \sup_{y\in S}\delta
(w,y)\leq \lambda \}$ is closed, and therefore compact since it is contained
in $\beta (y,\lambda )$\ for each $y\in S$. Consequently, there exists $x\in
A$ such that $\sup_{y\in S}\delta (x,y)=\inf_{w\in A}(\sup_{y\in S}\delta
(w,y))=\mathrm{Rad}(S)$.

Now suppose that there exists $x^{\prime }\neq x$ in $E$\ such that\ $%
\sup_{y\in S}\delta (x^{\prime },y)=\mathrm{Rad}(S)$, and consider $%
x^{\prime \prime }\in E$\ such that $\delta (x,x^{\prime \prime })=\delta
(x^{\prime \prime },x^{\prime })=\delta (x,x^{\prime })/2$. By (CVX$\lambda $%
), we have $\delta (x^{\prime \prime },y)<\sup (\delta (x,y),\delta
(x^{\prime },y))\leq \mathrm{Rad}(S)$ for each $y\in S$. It follows $%
\sup_{y\in S}\delta (x^{\prime \prime },y)\leq \mathrm{Rad}(S)$, and
therefore $\sup_{y\in S}\delta (x^{\prime \prime },y)=\mathrm{Rad}(S)$. As $%
S $ is compact, there exists $y\in S$\ such that $\delta (x^{\prime \prime
},y)=\mathrm{Rad}(S)$, which contradicts (CVX$\lambda $) since $\delta
(x,y)\leq \mathrm{Rad}(S)$ and $\delta (x^{\prime },y)\leq \mathrm{Rad}(S)$%
.~~$\blacksquare $\bigskip

\noindent \textbf{Proposition 1.10.} Suppose that $(E,\delta )$ is geodesic
and satisfies (CVX$\lambda $). Suppose that (FIX$\lambda \mu $) is true for
some $\mu \in 
\mathbb{R}
_{>0}$. Let $\Gamma =\{(\mathcal{C}_{1},C_{1}),...,(\mathcal{C}_{m},C_{m})\}$%
\ be a $1$-local rule defined with tiles of radius $<\lambda $. Consider $%
\xi \in 
\mathbb{R}
_{>0}$\ such that, for each $i\in \{1,...,m\}$, the union of the tiles of $%
\mathcal{C}_{i}$\ contains $\cup _{x\in C_{i}}\beta (x,\xi )$, and $p\in 
\mathbb{N}
^{\ast }$\ such that $p\xi \geq \mu $. Then no $\sigma \in G-\{\mathrm{Id}\}$
stabilizes the tiles of a $p$-configuration $(\mathcal{D},D)$\ which is
compatible with $\Gamma $. In particular, Theorem 1.3 and its corollaries
are true for each $(p+1)$-local rule $\Delta $\ which is compatible with $%
\Gamma $.\bigskip

\noindent \textbf{Proof.} Suppose that some $\sigma \in G-\{\mathrm{Id}\}$
stabilizes the tiles of a $p$-configuration $(\mathcal{D},D)$\ which is
compatible with $\Gamma $. As $(E,\delta )$ is geodesic and as each $(%
\mathcal{B}_{1}^{\mathcal{D}}(T),T)$ with $T\in \mathcal{B}_{p-1}^{\mathcal{D%
}}(D)$\ is isomorphic to some $(\mathcal{C}_{i},C_{i})$, we see by induction
on $0\leq k\leq p$ that $\beta (x,k\xi )$\ is contained in the union of the
tiles of $\mathcal{B}_{k}^{\mathcal{D}}(D)$. Consequently, $\beta (x,\mu
)\subset \beta (x,p\xi )$\ is contained in the union of the tiles of $%
\mathcal{D}$. For each $y\in \beta (x,\mu )$ and each $T\in \mathcal{D}$
such that $y\in T$, we have $\delta (y,x_{T})\leq \mathrm{Rad}(T)<\lambda $,
where $x_{T}$ is the unique point of $E$\ such that $\delta (z,x_{T})\leq 
\mathrm{Rad}(T)$\ for each $z\in T$.\ This contradicts (FIX$\lambda \mu $)
since $x_{T}\sigma =x_{T}$ for each $T\in \mathcal{D}$ by Lemma 1.9.~~$%
\blacksquare $\bigskip

\textbf{2. Local isomorphism and representation of relational structures by
tilings.}\bigskip

First we define local isomorphism and the extraction preorder $\Subset $ for
tilings and relational structures. The definitions for tilings are
classical. We note that, by the first part of Proposition 1.2, each subpatch
of a tiling is finite if and only if it is bounded.

We say that a tiling $\mathcal{T}$ satisfies the\emph{\ local isomorphism
property} if, for each finite subpatch $\mathcal{E}$ of $\mathcal{T}$, there
exists $k\in 
\mathbb{N}
^{\ast }$ such that each $\mathcal{B}_{k}^{\mathcal{T}}(T)$ contains a copy
of $\mathcal{E}$. Then, for each finite subpatch $\mathcal{E}$ of $\mathcal{T%
}$, there exists $\rho \in 
\mathbb{R}
_{>0}$\ such that each $\beta (x,\rho )\subset E$ contains a subpatch of $%
\mathcal{T}$ which is isomorphic to $\mathcal{E}$. The second part of
Proposition 1.2 implies that the converse is true if $E$ is weakly
homogeneous.

For any tilings $\mathcal{S},\mathcal{T}$, we write $\mathcal{S}\Subset 
\mathcal{T}$ if each finite subpatch of $\mathcal{S}$ is isomorphic to a
subpatch of $\mathcal{T}$. We say that $\mathcal{S}$\ and $\mathcal{T}$\ are 
\emph{locally isomorphic} if $\mathcal{S}\Subset \mathcal{T}$ and $\mathcal{T%
}\Subset \mathcal{S}$.

Now we give the definitions and the notations for relational structures,
which are similar to those in [7, pp. 107, 112, 113]. We consider a finite
relational language $\mathcal{L}$.

For each $\mathcal{L}$-structure $M$ and each $u\in M$, we define
inductively the subsets $B_{M}(u,h)$\ with $B_{M}(u,0)=\{u\}$ and, for $h\in 
\mathbb{N}
$,

\noindent $B_{M}(u,h+1)=B_{M}(u,h)\cup \{v\in M\mid $\ there exist $%
R(x_{1},...,x_{k})\in \mathcal{L}$, $u_{1},...,u_{k}\in M$ and $1\leq
i,j\leq k$\ such that $R(u_{1},...,u_{k})$, $u_{i}\in B_{M}(u,h)$ and $%
u_{j}=v\}$.

\noindent We say that $M$ is \emph{connected} if there exists $u\in M$ such
that $M=\cup _{h\in 
\mathbb{N}
}B_{M}(u,h)$.

We say that $M$ is \emph{locally finite} if $B_{M}(u,1)$ is finite for each $%
u\in M$. Then $B_{M}(u,h)$ is finite for any $u\in M$ and $h\in 
\mathbb{N}
$. We say that $M$ is \emph{uniformly locally finite} if there exists $r\in 
\mathbb{N}
^{\ast }$ such that $\left\vert B_{M}(u,1)\right\vert \leq r$ for each $u\in
M$. Then, for each $h\in 
\mathbb{N}
$, there exists $s\in 
\mathbb{N}
^{\ast }$ such that $\left\vert B_{M}(u,h)\right\vert \leq s$ for each $u\in
M$.\ We say that $M$ satisfies the \emph{local isomorphism property} if, for
any $u\in M$ and $h\in 
\mathbb{N}
$, there exists $k\in 
\mathbb{N}
$\ such that each $B_{M}(v,k)$ contains some $w$ with $(B_{M}(w,h),w)\cong
(B_{M}(u,h),u)$.

For any $\mathcal{L}$-structures $M,N$, we write $M\Subset N$ if, for any $%
u\in M$ and $h\in 
\mathbb{N}
$,

\noindent $\left\vert \left\{ v\in M\mid (B_{M}(v,h),v)\cong
(B_{M}(u,h),u)\right\} \right\vert \leq $

\noindent $\left\vert \left\{ w\in N\mid (B_{N}(w,h),w)\cong
(B_{M}(u,h),u)\right\} \right\vert $

\noindent or both sets are infinite. For $M,N$\ connected, we have $M\Subset
N$\ if and only if, for any $u\in M$ and $h\in 
\mathbb{N}
$, there exists $v\in N$\ such that $(B_{M}(u,h),u)\cong (B_{N}(v,h),v)$.

We say that $M$ and $N$ are \emph{locally isomorphic} if $M\Subset N$ and $%
N\Subset M$. By [7, Theorem 2.3], two locally finite $\mathcal{L}$%
-structures $M,N$ are elementarily equivalent if and only if they are
locally isomorphic.

Now we consider again the metric space with bounded closed subsets $%
(E,\delta )$, the group $G$ of bijective isometries of $E$, the integer $q$,
the $q$-local rule $\Delta =\{(\mathcal{D}_{1},D_{1}),...,(\mathcal{D}%
_{n},D_{n})\}$ such that each $\mathcal{B}_{q-1}^{\mathcal{D}_{i}}(D_{i})$
is fixed by no $\sigma \in G-\{\mathrm{Id}\}$, and the language $\mathcal{L}%
_{\Delta }$, which were introduced for Theorem 1.3. We call a $\Delta $\emph{%
-tiling} any tiling which satisfies $\Delta $.

We define a \emph{representation} of a $\mathcal{L}_{\Delta }$-structure $M$
as a pair $(\mathcal{E},f)$, where $\mathcal{E}$ is a $\Delta $\textbf{-}%
patch and $f:M\rightarrow \mathcal{E}$ is an isomorphism of $\mathcal{L}%
_{\Delta }$-structures. By Theorem 1.3, the representation is unique up to
isomorphism if it exists.

For each $\mathcal{L}_{\Delta }$-structure $M$ and each $u\in M$, we define
inductively the subsets $B_{h}^{M}(u)$ with $B_{0}^{M}(u)=\{u\}$ and, for
each $h\in 
\mathbb{N}
$,

\noindent $B_{h+1}^{M}(u)=B_{h}^{M}(u)\cup \{v\in M\mid $\ there exist $%
1\leq i\leq n$, $v_{1},...,v_{p(i)}\in M$ and $1\leq j,k\leq p(i)$\ such
that $R_{i}(v_{1},...,v_{p(i)})$, $v_{j}\in B_{h}^{M}(u)$, $v_{k}=v$\ and $%
D_{i,j}\cap D_{i,k}\neq \emptyset \}$

\noindent where the sets $D_{i,j}$ are the tiles used for the definition of $%
\mathcal{L}_{\Delta }$.

For each set $\mathcal{E}$ of tiles, each $T\in \mathcal{E}$ and each $h\in 
\mathbb{N}
$, the $\mathcal{L}_{\Delta }$-structure $M$ defined on $\mathcal{E}$
satisfies $B_{h}^{M}(T)\subset \mathcal{B}_{h}^{\mathcal{E}}(T)$, and $%
B_{h}^{M}(T)=\mathcal{B}_{h}^{\mathcal{E}}(T)$ if $\mathcal{E}$ is a $\Delta 
$-tiling. The notation $B_{h}^{M}(u)$ should not be confused with $%
B_{M}(u,h) $. For each $\mathcal{L}_{\Delta }$-structure $M$, each $u\in M$
and each $h\in 
\mathbb{N}
$, we have $B_{h}^{M}(u)\subset B_{M}(u,h)\subset B_{2qh}^{M}(u)$, and $%
B_{M}(u,h)=B_{2qh}^{M}(u)$ if $M$ can be represented by a $\Delta $-tiling.

Now we introduce something analogous to the notion of local rule considered
for tilings. A \emph{local rule} for a finite relational language $\mathcal{L%
}$ is defined by specifying an integer $r\in 
\mathbb{N}
^{\ast }$ and a finite sequence of pairs $(M_{1},x_{1}),...,(M_{k},x_{k})$
with $M_{1},...,M_{k}$ finite $\mathcal{L}$-structures and $%
M_{i}=B_{M_{i}}(x_{i},r)$ for $1\leq i\leq k$. We say that a $\mathcal{L}$%
-structure $M$ $\emph{satisfies}$ that rule if each $(B_{M}(x,r),x)$ is
isomorphic to one of the pairs $(M_{i},x_{i})$.

Any such rule can be expressed by a first-order sentence. For each $s\in 
\mathbb{N}
^{\ast }$, the pairs $(B_{M}(x,s),x)$ for $M$ satisfying that rule and $x\in
M$\ are finite and\ fall in finitely many isomorphism classes. In
particular, any $\mathcal{L}$-structure which satisfies a local rule\ is
uniformly locally finite.

Theorem 2.1 below implies that there exists a generally infinite set of
local rules which characterizes, among the connected $\mathcal{L}_{\Delta }$%
-structures, those which can be represented by $\Delta $-tilings. The
problem of the characterization of representable $\mathcal{L}_{\Delta }$%
-structures by a finite set of local rules will be considered with Theorem
2.7 and the examples in Section 4.

The following property of a $\mathcal{L}_{\Delta }$-structure $M$ is true,
in particular, if $M$ can be represented by a $\Delta $\textbf{-}tiling:

\noindent (P) For each $u\in M$,\ there exist $1\leq i\leq n$\ and $%
u_{2},...,u_{p(i)}\in M$ such that $R_{i}(u,u_{2},...,u_{p(i)})$.\bigskip

\noindent \textbf{Theorem 2.1.} For each connected $\mathcal{L}_{\Delta }$%
-structure $M$, the following properties are equivalent:

\noindent 1) $M$\ can be represented by a $\Delta $-tiling;

\noindent 2) For each $u\in M$ and each $r\in 
\mathbb{N}
^{\ast }$, we have $(B_{M}(u,r),u)\cong (B_{\mathcal{T}}(T,r),T)$ for a tile 
$T$ of a $\Delta $-tiling $\mathcal{T}$;

\noindent 3) $M$\ satisfies (P) and, for each $u\in M$\ and each integer $%
h\geq q$, there exists a representation of $B_{h}^{M}(u)$ by a $\Delta $%
-patch.\bigskip

\noindent \textbf{Proof.} The property 1) clearly implies 2).

Suppose that 2) is true. Consider any element $u\in M$\ and any integer $%
h\geq q$. Then, for each $\Delta $-tiling $\mathcal{T}$ and each $T\in 
\mathcal{T}$, each isomorphism from $(B_{M}(u,h),u)\ $to $(B_{\mathcal{T}%
}(T,h),T)$ induces an isomorphism from $(B_{h}^{M}(u),u)\subset
(B_{M}(u,h),u)$ to $(B_{h}^{\mathcal{T}}(T),T)$.\ Such an isomorphism gives
a representation of $B_{h}^{M}(u)$ since $B_{h}^{\mathcal{T}}(T)=\mathcal{B}%
_{h}^{\mathcal{T}}(T)$\ is a $\Delta $-patch.\ It also induces an
isomorphism from $(B_{q}^{M}(u),u)$ to $(B_{q}^{\mathcal{T}}(T),T)$, which
implies that\ there exist $1\leq i\leq n$\ and $u_{2},...,u_{p(i)}\in M$
such that $R_{i}(u,u_{2},...,u_{p(i)})$.

Now suppose that 3) is true and fix $u\in M$. For each integer $h\geq q$,
consider a representation $(\mathcal{E}_{h},f_{h})$ of $B_{h}^{M}(u)$ by a $%
\Delta $-patch. For $k\geq h\geq q$, as $f_{h}^{-1}f_{k}$ is a homomorphism
of $\mathcal{L}_{\Delta }$-structures, Theorem 1.3 implies that there exists
a unique $\sigma _{h,k}\in G$ such that $Sf_{h}^{-1}f_{k}=S\sigma _{h,k}$
for each $S\in \mathcal{E}_{h}$. So we can suppose that, for $k\geq h\geq q$%
, we have $\mathcal{E}_{h}\subset \mathcal{E}_{k}$ and $f_{h}$ is the
restriction of $f_{k}$\ to $B_{h}^{M}(u)$. Then $f=\cup _{h\geq q}f_{h}$ is
a $\mathcal{L}_{\Delta }$\textbf{-}isomorphism from $M$ to $\mathcal{E}=\cup
_{h\geq q}\mathcal{E}_{h}$. As $M$ satisfies (P), it follows that $\mathcal{E%
}$\ satisfies $\Delta $, and $\mathcal{E}$\ is a tiling by Corollary 1.7.~~$%
\blacksquare $\bigskip

\noindent \textbf{Corollary 2.2.} For any connected $\mathcal{L}_{\Delta }$%
-structures $M\Subset N$, if $N$ can be represented by a $\Delta $-tiling,
then $M$ can also be represented by a $\Delta $-tiling.\bigskip

Similarly to [7], we have:\bigskip

\noindent \textbf{Proposition 2.3.} 1) Any $\Delta $-tiling satisfies the\
local isomorphism property as a tiling if and only if it satisfies the\
local isomorphism property as a $\mathcal{L}_{\Delta }$-structure.

\noindent 2) Two $\Delta $-tilings $\mathcal{S},\mathcal{T}$ satisfy $%
\mathcal{S}\Subset \mathcal{T}$ (resp. are locally isomorphic) as tilings if
and only if they satisfy $\mathcal{S}\Subset \mathcal{T}$ (resp. are locally
isomorphic) as $\mathcal{L}_{\Delta }$-structures.\bigskip

\noindent \textbf{Proof.} The following facts will be used in the proofs of
1) and 2):

\noindent a) We have $B_{\mathcal{T}}(T,h)=B_{2qh}^{\mathcal{T}}(T)=\mathcal{%
B}_{2qh}^{\mathcal{T}}(T)$ for each $\Delta $\textbf{-}tiling $\mathcal{T}$,
each $T\in \mathcal{T}$ and each $h\in 
\mathbb{N}
$.

\noindent b) Theorem 1.3 implies that, for any $\Delta $-tilings $\mathcal{S}%
,\mathcal{T}$, each $S\in \mathcal{S}$, each $T\in \mathcal{T}$ and each
integer $k\geq q$, the $\Delta $-patches\ $\mathcal{B}_{k}^{\mathcal{S}}(S)$
and $\mathcal{B}_{k}^{\mathcal{T}}(T)$ are isomorphic as sets of tiles if
and only if they are isomorphic as $\mathcal{L}_{\Delta }$-structures.

\noindent c) For any $\Delta $-tilings $\mathcal{S},\mathcal{T}$, each $S\in 
\mathcal{S}$, each $k\in 
\mathbb{N}
^{\ast }$ and each $\sigma \in G$ such that\ $\mathcal{B}_{k}^{\mathcal{S}%
}(S)\sigma \subset \mathcal{T}$, we have\ $\mathcal{B}_{k}^{\mathcal{S}%
}(S)\sigma =\mathcal{B}_{k}^{\mathcal{T}}(S\sigma )$.

\noindent d) Each finite subpatch of a $\Delta $\textbf{-}tiling $\mathcal{T}
$ is contained in some $\mathcal{B}_{h}^{\mathcal{T}}(T)$.

First we prove 1). The facts a), b), c) above imply that $\mathcal{T}$
satisfies the\ local isomorphism property as a $\mathcal{L}_{\Delta }$%
-structure if and only if, for each $h\in 
\mathbb{N}
^{\ast }$ and each $S\in \mathcal{T}$,\ there exists $k\in 
\mathbb{N}
^{\ast }$ such that each $\mathcal{B}_{k}^{\mathcal{T}}(T)$ contains a copy
of $\mathcal{B}_{h}^{\mathcal{T}}(S)$. By d), the last property is true if
and only if $\mathcal{T}$ satisfies the\ local isomorphism property as a
tiling.

Now we prove 2). We only show the first statement, since the second one is
an immediate consequence. The facts a), b), c) above imply that $\mathcal{S}$
and $\mathcal{T}$ satisfy $\mathcal{S}\Subset \mathcal{T}$ as $\mathcal{L}%
_{\Delta }$-structures if and only if $\mathcal{T}$ contains a copy of $%
\mathcal{B}_{h}^{\mathcal{S}}(S)$ for each $h\in 
\mathbb{N}
^{\ast }$ and each $S\in \mathcal{S}$. By d), the last property is true if
and only if $\mathcal{S}$ and $\mathcal{T}$ satisfy $\mathcal{S}\Subset 
\mathcal{T}$ as tilings.~~$\blacksquare $\bigskip

By Theorem 2.1, for each sequence $(B_{M_{i}}(u_{i},r_{i}),u_{i})_{i\in 
\mathbb{N}
}=(B_{2qr_{i}}^{M_{i}}(u_{i}),u_{i})_{i\in 
\mathbb{N}
}$ of pairs taken in $\mathcal{L}_{\Delta }$-structures associated to $%
\Delta $\textbf{-}tilings, with $r_{i}<r_{j}$ for $i<j$, the inductive limit
relative to any sequence of isomorphisms

\noindent $f_{i}:(B_{M_{i}}(u_{i},r_{i}),u_{i})\rightarrow
(B_{M_{i+1}}(u_{i+1},r_{i}),u_{i+1})\subset
(B_{M_{i+1}}(u_{i+1},r_{i+1}),u_{i+1})$

\noindent is a pair $(M,x)$\ with $M$\ a $\mathcal{L}_{\Delta }$-structure
associated to a $\Delta $\textbf{-}tiling. Using this fact, together with
Corollary 2.2 and Proposition 2.3, we see that Corollary 2.4 (resp. 2.5,
resp. 2.6) below, in the same way as [7, Corollary 3.5] (resp. [7, Corollary
3.6], resp. [7, Corollary 3.7]) is a consequence of [7, Corollary 3.2]
(resp. [7, Proposition 3.3], resp. [7, Proposition 3.4]) and its
proof.\bigskip

\noindent \textbf{Corollary 2.4.} Any $\Delta $\textbf{-}tiling is minimal
for $\Subset $ if and only if it satisfies the local isomorphism
property.\bigskip

\noindent \textbf{Corollary 2.5.} For each $\Delta $\textbf{-}tiling $%
\mathcal{S}$, there exists a $\Delta $\textbf{-}tiling $\mathcal{T}\Subset 
\mathcal{S}$ which is minimal for $\Subset $.\bigskip

\noindent \textbf{Corollary 2.6.} For each $\Delta $\textbf{-}tiling $%
\mathcal{S}$, we have:

\noindent 1) If there are finitely many equivalence classes of elements of $%
\mathcal{S}$ modulo the isometries $\sigma \in G$ such that $\mathcal{S}%
\sigma =\mathcal{S}$, then any $\Delta $\textbf{-}tiling $\mathcal{T}\Subset 
\mathcal{S}$\ is isomorphic to $\mathcal{S}$.

\noindent 2) If $\mathcal{S}$ is minimal for $\Subset $, and if there are
infinitely many equivalence classes of elements of $\mathcal{S}$ modulo the
isometries $\sigma \in G$ such that $\mathcal{S}\sigma =\mathcal{S}$, then
there exist $2^{\omega }$ pairwise nonisomorphic $\Delta $\textbf{-}tilings
which are locally isomorphic to $\mathcal{S}$.\bigskip

\noindent \textbf{Remark.}\ Corollary 2.6 implies [11, Theorem, p. 356]. In
[11], Radin and Wolff considered tilings of the euclidean spaces $%
\mathbb{R}
^{n}$, and isomorphism was defined modulo an arbitrary group of
isometries.\bigskip

Theorem 2.7 below is analogous to [7, Theorem 5.2]. Here, we suppose that $%
(E,\delta )$ is weakly homogeneous and that, for any $x,y\in E$ and each $%
\eta \in 
\mathbb{R}
_{>0}$, $\beta (x,\delta (x,y))\cap \beta (y,\eta )$\ is connected and
contains a point $z$ with $\delta (x,z)<\delta (x,y)$.

These conditions are true for each geodesic space which satisfies the
properties (CVX$\nu $) considered at the end of Section 1. In each $%
\mathbb{R}
^{k}$, they are also true for the distances defined with

\noindent $\delta ((x_{1},...,x_{k}),(y_{1},...,y_{k}))=\sup (\left\vert
y_{1}-x_{1}\right\vert ,...,\left\vert y_{k}-x_{k}\right\vert )$\ or

\noindent $\delta ((x_{1},...,x_{k}),(y_{1},...,y_{k}))=\left\vert
y_{1}-x_{1}\right\vert +...+\left\vert y_{k}-x_{k}\right\vert )$.

On the other hand, the connectedness condition is not true if $E$ is the
surface of a cylinder of infinite length and if $\delta $ is defined by
considering $E$ as a quotient of the euclidean space $%
\mathbb{R}
^{2}$. Actually, Theorem 2.7 is not true in that case since, for some local
rules $\Delta $, it is possible to find a $\mathcal{L}_{\Delta }$-structure $%
M$ which can only be represented by a tiling of $%
\mathbb{R}
^{2}$, while the substructures $B_{r}^{M}(u)$ considered in Theorem 2.7 are
small enough to be represented in $(E,\delta )$.

Again, we consider $\mathcal{L}_{\Delta }$-structures which satisfy the
property (P) introduced for Theorem 2.1. We denote by $\lambda $ the maximum
of the diameters of the tiles in $\Delta $, and $\xi $ the largest real
number such that, for each $i\in \{1,...,n\}$, the union of the tiles of $%
\mathcal{B}_{1}^{\mathcal{D}_{i}}(D_{i})$\ contains $\cup _{x\in D_{i}}\beta
(x,\xi )$. For each $r\in 
\mathbb{N}
$, we write $\omega (\xi ,r)=\inf_{x\in E}\omega (x,\xi ,r)$, as in the
definition of weakly homogeneous spaces.\bigskip

\noindent \textbf{Theorem 2.7.} Suppose that $(E,\delta )$ is weakly
homogeneous and that, for any $x,y\in E$ and each $\eta \in 
\mathbb{R}
_{>0}$, $\beta (x,\delta (x,y))\cap \beta (y,\eta )$\ is connected and
contains a point $z$ with $\delta (x,z)<\delta (x,y)$. Define (P), $\lambda $%
\ and $\xi $\ as above. Consider an integer $r\geq q+1$\ such that $\omega
(\xi ,r-q-1)\geq (4q+1)\lambda $. Let $M$ be a connected $\mathcal{L}%
_{\Delta }$-structure which satisfies (P) and such that each $B_{r}^{M}(u)$
can be represented by a $\Delta $\textbf{-}patch. Then there exist a $\Delta 
$\textbf{-}tiling $\mathcal{T}$ and a surjective $\mathcal{L}_{\Delta }$%
-homomorphism $\varphi :\mathcal{T}\rightarrow M$ which induces a $\mathcal{L%
}_{\Delta }$-isomorphism from $\mathcal{B}_{r}^{\mathcal{T}}(T)$ to $%
B_{r}^{M}(T\varphi )$ for each\ $T\in \mathcal{T}$.\bigskip

\noindent \textbf{Remark.}\ The pair $(\mathcal{T},\varphi )$ satisfies $%
S\varphi \neq T\varphi $ for each $S\in \mathcal{T}$ and each\ $T\in 
\mathcal{B}_{2r}^{\mathcal{T}}(S)-\{S\}$, since there exists $U\in \mathcal{T%
}$\ such that\ $S,T\in \mathcal{B}_{r}^{\mathcal{T}}(U)$.\bigskip

\noindent \textbf{Remark.}\ If $(E,\delta )$ is geodesic, we can take any
integer $r\geq q+1+(4q+1)\lambda /\xi $.\bigskip

\noindent \textbf{Remark.} It follows that there exists a $\Delta $\textbf{-}%
tiling if there exists a $\mathcal{L}_{\Delta }$-structure $M$ which
satisfies (P) and such that each $B_{r}^{M}(u)$ can be represented by a $%
\Delta $\textbf{-}patch.\bigskip

\noindent \textbf{Lemma 2.7.1.} For each\ $u\in M$, each representation $(%
\mathcal{E},f)$\ of $B_{r}^{M}(u)$ and each $x\in uf$, any\ $T\in \mathcal{E}
$ such that $T\cap \beta (x,(4q+1)\lambda )\neq \emptyset $ belongs to $%
\mathrm{I}_{\Delta }(\mathcal{E})$.\bigskip

\noindent \textbf{Proof of Lemma 2.7.1.} As $M$\ satisfies (P), for each $%
v\in B_{r-q}^{M}(u)$, we have $vf\in \mathrm{I}_{\Delta }(\mathcal{E})$ and $%
\cup _{y\in vf}\beta (y,\xi )$ is contained in the union of the tiles of $%
B_{1}^{M}(v)f$. An induction on $k$ shows that $\beta (x,\omega (\xi ,k))$
is contained in the union of the tiles of $B_{k}^{M}(u)f$ for $0\leq k\leq
r-q$. In particular, $\beta (x,(4q+1)\lambda )\subset \beta (x,\omega (\xi
,r-q-1))$ is contained in the union of the tiles of $B_{r-q-1}^{M}(u)f$.
Consequently, any\ $T\in \mathcal{E}$ such that $T\cap \beta
(x,(4q+1)\lambda )\neq \emptyset $ belongs to $B_{r-q}^{M}(u)f$, and
therefore belongs to $\mathrm{I}_{\Delta }(\mathcal{E})$.~~$\blacksquare $%
\bigskip

\noindent \textbf{Proof of Theorem 2.7.} We consider an element\ $u\in M$, a
representation $(\mathcal{G},g)$\ of $B_{r}^{M}(u)$ and a point $x\in ug$.
We denote by $\Omega $ the set of all pairs $(\mathcal{E},\varphi )$\ such
that:

\noindent a) $\mathcal{E}$\ is a finite $\Delta $\textbf{-}patch, $x$
belongs to a tile of $\mathcal{E}$, any\ $T\in \mathcal{E}$ such that $x\in
T $ belongs to $\mathrm{I}_{\Delta }(\mathcal{E})$, and $\varphi $ is a $%
\mathcal{L}_{\Delta }$-homomorphism from $\mathcal{E}$\ to $M$;

\noindent b) for each $S\in \mathcal{E}$, there exists $T\in \mathrm{I}%
_{\Delta }(\mathcal{E})$ such that $S\in \mathcal{B}_{q}^{\mathcal{E}}(T)$
and $T\cap \beta (x,\rho _{\mathcal{E}})\neq \emptyset $, where $\rho _{%
\mathcal{E}}=\sup \{\rho \in 
\mathbb{R}
_{\geq 0}\mid $ any\ $U\in \mathcal{E}$ such that $U\cap \beta (x,\rho )\neq
\emptyset $ belongs to $\mathrm{I}_{\Delta }(\mathcal{E})\}$.

For each $(\mathcal{E},\varphi )\in \Omega $, $\beta (x,\rho _{\mathcal{E}})$
is covered by the tiles of $\mathcal{E}$: For each $y\in E$ such that $%
\delta (x,y)=\rho _{\mathcal{E}}$, the second condition on $(E,\delta )$
implies that $y$\ is the limit of a sequence $(y_{k})_{k\in 
\mathbb{N}
}\subset \cup _{\zeta <\rho _{\mathcal{E}}}\beta (x,\zeta )$. As $\mathrm{I}%
_{\Delta }(\mathcal{E})$\ is finite, infinitely many $y_{k}$\ belong to the
same $T\in \mathrm{I}_{\Delta }(\mathcal{E})$\ and $T$ also contains $y$.

On the other hand, there exists $y\in E$ with $\delta (x,y)=\rho _{\mathcal{E%
}}$ which belongs to some $S\in \mathcal{E}-\mathrm{I}_{\Delta }(\mathcal{E}%
) $: For any sequences $(y_{k})_{k\in 
\mathbb{N}
}\subset E$\ and $(S_{k})_{k\in 
\mathbb{N}
}\subset \mathcal{E}-\mathrm{I}_{\Delta }(\mathcal{E})$ such that $y_{k}\in
S_{k}$\ for each $k\in 
\mathbb{N}
$ and $\lim_{k\rightarrow +\infty }\delta (y_{k},\beta (x,\rho _{\mathcal{E}%
}))=0$,\ there exists a subsequence of $(y_{k})_{k\in 
\mathbb{N}
}$ which converges to a point $y\in \beta (x,\rho _{\mathcal{E}})$. As $%
\mathcal{E}-\mathrm{I}_{\Delta }(\mathcal{E})$\ is finite, infinitely many $%
y_{k}$\ in this subsequence belong to the same $S_{h}$\ and $S_{h}$ also
contains $y$.

By Lemma 2.7.1, $\Omega $ contains $(\mathcal{E}_{0},\varphi _{0})$, where $%
\mathcal{E}_{0}$ is the union of the subsets $\mathcal{B}_{q}^{\mathcal{G}%
}(T)$ for $T\in \mathcal{G}$\ such that $T\cap \beta (x,(4q+1)\lambda )\neq
\emptyset $, and $\varphi _{0}$ is the restriction of $g^{-1}$ to $\mathcal{E%
}_{0}$. We have $\rho _{\mathcal{E}_{0}}>(4q+1)\lambda $.

For any $(\mathcal{E},\varphi ),(\mathcal{F},\psi )\in \Omega $, we write $(%
\mathcal{E},\varphi )\leq (\mathcal{F},\psi )$ if $\mathcal{E}\subset 
\mathcal{F}$\ and if $\varphi $\ is the restriction of $\psi $ to $\mathcal{E%
}$. It suffices to prove the two following claims:

1.\textit{\ The conclusion of Theorem 2.7 is true if }$\Omega $\textit{\
contains some strictly increasing }$(\mathcal{E}_{i},\varphi _{i})_{i\in 
\mathbb{N}
}$.

First we observe that, for each $i\in 
\mathbb{N}
$, there exists $T\in \mathrm{I}_{\Delta }(\mathcal{E}_{i+1})-\mathrm{I}%
_{\Delta }(\mathcal{E}_{i})$ such that $T\cap \beta (x,\rho _{\mathcal{E}%
_{i}})\neq \emptyset $: If $\rho _{\mathcal{E}_{i+1}}>\rho _{\mathcal{E}%
_{i}} $,\ then any $T\in \mathcal{E}_{i}-\mathrm{I}_{\Delta }(\mathcal{E}%
_{i})$ such that $T\cap \beta (x,\rho _{\mathcal{E}_{i}})\neq \emptyset $\
belongs to$\ \mathrm{I}_{\Delta }(\mathcal{E}_{i+1})$. If $\rho _{\mathcal{E}%
_{i+1}}=\rho _{\mathcal{E}_{i}}$, then, for each $S\in \mathcal{E}_{i+1}-%
\mathcal{E}_{i}$,$\ $there\ exists $T\in \mathrm{I}_{\Delta }(\mathcal{E}%
_{i+1})\ $such that $S\in \mathcal{B}_{q}^{\mathcal{E}_{i+1}}(T)$ and $T\cap
\beta (x,\rho _{\mathcal{E}_{i}})=\emptyset $; any such $T$ cannot belong to 
$\mathrm{I}_{\Delta }(\mathcal{E}_{i})$ since $\mathcal{B}_{q}^{\mathcal{E}%
_{i+1}}(T)-\mathcal{B}_{q}^{\mathcal{E}_{i}}(T)\neq \emptyset $.

Then we use the function $\omega (\xi ,n)$ in order to prove that $\rho _{%
\mathcal{E}_{i}}$\ tends to infinity. As $(E,\delta )$ is weakly
homogeneous, $\omega (\xi ,n)$\ tends to infinity with $n$. Consequently, it
suffices to show that, for any $i,n\in 
\mathbb{N}
$ such that $\rho _{\mathcal{E}_{i}}>\omega (\xi ,n)$, there exists $j>i$
such that $\rho _{\mathcal{E}_{j}}>\omega (\xi ,n+1)$.

For each $j\geq i$ such that $\rho _{\mathcal{E}_{j}}\leq \omega (\xi ,n+1)$%
, we consider$\ T_{j}\in \mathrm{I}_{\Delta }(\mathcal{E}_{j+1})-\mathrm{I}%
_{\Delta }(\mathcal{E}_{j})$ such that $T_{j}\cap \beta (x,\rho _{\mathcal{E}%
_{j}})\neq \emptyset $, and $y_{j}\in T_{j}\cap \beta (x,\rho _{\mathcal{E}%
_{j}})$. There exist $z_{j}\in \beta (x,\omega (\xi ,n))$ such that $\delta
(y_{j},z_{j})\leq \xi $, and $U_{j}\in \mathcal{B}_{1}^{\mathcal{E}j}(T_{j})$
such that $z_{j}\in U_{j}$. We have $U_{j}\in \mathrm{I}_{\Delta }(\mathcal{E%
}_{i})$ since $z_{j}\in \beta (x,\omega (\xi ,n))$ and\ $\omega (\xi
,n)<\rho _{\mathcal{E}_{i}}$. It follows $T_{j}\in \mathcal{E}_{i}$. As $%
\mathcal{E}_{i}$ is finite, there exist at most finitely\ many such $T_{j}$,
and therefore finitely\ many integers $j>i$ such that $\rho _{\mathcal{E}%
_{j}}\leq \omega (\xi ,n+1)$.

As $\rho _{\mathcal{E}_{i}}$\ tends to infinity,\ $\mathcal{T}=\cup _{i\in 
\mathbb{N}
}\mathcal{E}_{i}$ satisfies $\Delta $, and $\mathcal{T}$ is a $\Delta $%
-tiling by Corollary 1.7. The inductive limit $\varphi $ of the maps $%
\varphi _{i}$\ is a $\mathcal{L}_{\Delta }$-homomorphism from $\mathcal{T}$
to $M$. In particular, we have\ $\mathcal{B}_{k}^{\mathcal{T}}(T)\varphi
\subset B_{k}^{M}(T\varphi )$ for each\ $T\in \mathcal{T}$ and each $k\in 
\mathbb{N}
$.

Now we show that, for each\ $T\in \mathcal{T}$, $\varphi $ induces a $%
\mathcal{L}_{\Delta }$-isomorphism from $\mathcal{B}_{r}^{\mathcal{T}}(T)$
to $B_{r}^{M}(T\varphi )$; it follows that $\varphi $\ is surjective since $%
M $ is connected.

We consider a representation $(\mathcal{H},h)$\ of $B_{r}^{M}(T\varphi )$.
As $\mathcal{B}_{r}^{\mathcal{T}}(T)$ and $\mathcal{H}$ are $\Delta $%
-patches, Theorem 1.3 implies that there exists $\sigma \in G$ such that $%
S\sigma =S\varphi h$ for each $S\in \mathcal{B}_{r}^{\mathcal{T}}(T)$. As $%
\mathcal{T}$ is a $\Delta $-tiling and $\mathcal{H}$ is a $\Delta $-patch, $%
\sigma $\ induces an isomorphism from $\mathcal{B}_{q}^{\mathcal{T}}(U)$ to $%
\mathcal{B}_{q}^{\mathcal{H}}(U\sigma )$ for each $U\in \mathcal{B}_{r-q}^{%
\mathcal{T}}(T)\subset I_{\Delta }(\mathcal{B}_{r}^{\mathcal{T}}(T))$.
Consequently, we have $\mathcal{B}_{r}^{\mathcal{T}}(T)\sigma =\mathcal{B}%
_{r}^{\mathcal{H}}(T\sigma )=\mathcal{B}_{r}^{\mathcal{H}}(T\varphi
h)=B_{r}^{M}(T\varphi )h=\mathcal{H}$. It follows that $\sigma $ induces an
isomorphism from $\mathcal{B}_{r}^{\mathcal{T}}(T)$ to $\mathcal{H}$ and $%
\varphi $ induces a $\mathcal{L}_{\Delta }$-isomorphism from $\mathcal{B}%
_{r}^{\mathcal{T}}(T)$ to $B_{r}^{M}(T\varphi )$.

2. \textit{For each }$(\mathcal{E},\varphi )\in \Omega $\textit{, there
exists }$(\mathcal{F},\psi )>(\mathcal{E},\varphi )$\textit{\ in }$\Omega $.

We consider $y\in E$ with $\delta (x,y)=\rho _{\mathcal{E}}$ which belongs
to some $S\in \mathcal{E}-\mathrm{I}_{\Delta }(\mathcal{E})$, and $T\in 
\mathrm{I}_{\Delta }(\mathcal{E})\ $such that $y\in T$. There exist a
representation of $B_{r}^{M}(T\varphi )$, and therefore a $\Delta $-patch $%
\mathcal{E}^{\prime }$ and a $\mathcal{L}_{\Delta }$-isomorphism $\varphi
^{\prime }$ from $\mathcal{E}^{\prime }$\ to $B_{r}^{M}(T\varphi )$. We have 
$\mathcal{B}_{q}^{\mathcal{E}}(T)\varphi \subset B_{q}^{M}(T\varphi )$ since 
$T$ belongs to $\mathrm{I}_{\Delta }(\mathcal{E})$, and $B_{q}^{M}(T\varphi
)\varphi ^{\prime -1}\subset \mathcal{B}_{q}^{\mathcal{E}^{\prime
}}(T\varphi \varphi ^{\prime -1})$\ since $\varphi ^{\prime }$ is a $%
\mathcal{L}_{\Delta }$-isomorphism. The map\ $\mathcal{B}_{q}^{\mathcal{E}%
}(T)\rightarrow \mathcal{B}_{q}^{\mathcal{E}^{\prime }}(T\varphi \varphi
^{\prime -1}):U\rightarrow U\varphi \varphi ^{\prime -1}$ is a $\mathcal{L}%
_{\Delta }$-homomorphism. As $\mathcal{B}_{q}^{\mathcal{E}}(T)$ and $%
\mathcal{B}_{q}^{\mathcal{E}^{\prime }}(T\varphi \varphi ^{\prime -1})$ are $%
\Delta $-patches, Theorem 1.3 implies that there exists $\sigma \in G$ such
that $U\sigma =U\varphi \varphi ^{\prime -1}$ for each $U\in \mathcal{B}%
_{q}^{\mathcal{E}}(T)$. So we can suppose for the remainder of the proof
that $\mathcal{B}_{q}^{\mathcal{E}}(T)\subset \mathcal{E}^{\prime }$ and $%
U\varphi =U\varphi ^{\prime }\ $for each $U\in \mathcal{B}_{q}^{\mathcal{E}%
}(T)$.

We denote by $\mathcal{F}\ $the union of $\mathcal{E}\ $and the subsets $%
\mathcal{B}_{q}^{\mathcal{E}^{\prime }}(U)$ for $U\in \mathcal{E}^{\prime }$
such that $y\in U$. We consider the map $\psi :\mathcal{F}\rightarrow M$
with $U\psi =U\varphi $\ for $U\in \mathcal{E}$ and $U\psi =U\varphi
^{\prime }$\ for $U\in \mathcal{F}-\mathcal{E}$.

First we prove that, for each $U\in \mathcal{F}$, we have $\mathcal{B}_{q}^{%
\mathcal{F}}(U)\subset \mathcal{E}$ and $V\psi =V\varphi $ for each $V\in 
\mathcal{B}_{q}^{\mathcal{F}}(U)$, or $\mathcal{B}_{q}^{\mathcal{F}%
}(U)\subset \mathcal{E}^{\prime }$ and $V\psi =V\varphi ^{\prime }$ for each 
$V\in \mathcal{B}_{q}^{\mathcal{F}}(U)$.

If $U\cap \beta (y,2q\lambda )=\emptyset $, then the tiles in $\mathcal{B}%
_{q}^{\mathcal{F}}(U)$ cannot belong to $\mathcal{F}-\mathcal{E}$ since they
contain no point in $\beta (y,q\lambda )$.\ Consequently, we have $\mathcal{B%
}_{q}^{\mathcal{F}}(U)\subset \mathcal{E}$ and $V\psi =V\varphi $ for each $%
V\in \mathcal{B}_{q}^{\mathcal{F}}(U)$.

If $U\cap \beta (y,2q\lambda )\neq \emptyset $, then we prove that $\mathcal{%
B}_{q}^{\mathcal{F}}(U)\subset \mathcal{E}^{\prime }$ and $V\psi =V\varphi
^{\prime }$ for each $V\in \mathcal{B}_{q}^{\mathcal{F}}(U)$. For each $V\in 
\mathcal{B}_{q}^{\mathcal{F}}(U)$, we have $V\cap \beta (y,3q\lambda )\neq
\emptyset $. So it suffices to show that each $V\in \mathcal{E}$ with $V\cap
\beta (y,3q\lambda )\neq \emptyset $ satisfies $V\in \mathcal{E}^{\prime }$
and $V\varphi =V\varphi ^{\prime }$. We consider $W\in \mathrm{I}_{\Delta }(%
\mathcal{E})$\ such that $V\in \mathcal{B}_{q}^{\mathcal{E}}(W)$ and $W\cap
\beta (x,\rho _{\mathcal{E}})\neq \emptyset $. As $W\cap \beta (y,4q\lambda
)\neq \emptyset $,\ the properties $V\in \mathcal{E}^{\prime }$ and $%
V\varphi =V\varphi ^{\prime }$\ follow from Lemma 2.7.2 below.

Consequently, $\psi $ is a $\mathcal{L}_{\Delta }$-homomorphism and, for
each $U\in \mathcal{F}$, there exists $i\in \{1,...,n\}$\ such that $(%
\mathcal{B}_{q}^{\mathcal{F}}(U),U)\leq (\mathcal{D}_{i},D_{i})$. Moreover,
the definition of $(\mathcal{F},\psi )$\ implies that, for each $V\in 
\mathcal{F}$, there exists $U\in \mathrm{I}_{\Delta }(\mathcal{E})$ such
that $V\in \mathcal{B}_{q}^{\mathcal{E}}(U)$ and $U\cap \beta (x,\rho _{%
\mathcal{E}})\neq \emptyset $, or $U\in \mathrm{I}_{\Delta }(\mathcal{E}%
^{\prime })$ such that $V\in \mathcal{B}_{q}^{\mathcal{E}^{\prime }}(U)$ and 
$y\in U$; in both cases, we have $U\in \mathrm{I}_{\Delta }(\mathcal{F})$, $%
U\cap \beta (x,\rho _{\mathcal{E}})\neq \emptyset $ and $V\in \mathcal{B}%
_{q}^{\mathcal{F}}(U)$. Consequently, $\mathcal{F}$\ is a finite $\Delta $%
\textbf{-}patch and satisfies the conditions of the definition of $\Omega $.$%
~~\blacksquare $\bigskip

\noindent \textbf{Lemma 2.7.2.} For each $W\in \mathrm{I}_{\Delta }(\mathcal{%
E})$\ such that $W\cap \beta (x,\rho _{\mathcal{E}})\neq \emptyset $ and $%
W\cap \beta (y,4q\lambda )\neq \emptyset $,\ we have $\mathcal{B}_{q}^{%
\mathcal{E}}(W)\subset \mathcal{E}^{\prime }$ and $X\varphi =X\varphi
^{\prime }$\ for each $X\in \mathcal{B}_{q}^{\mathcal{E}}(W)$.\bigskip

\noindent \textbf{Proof of Lemma 2.7.2.} We fix $z\in W\cap \beta (x,\rho _{%
\mathcal{E}})$. We have $z\in \beta (y,(4q+1)\lambda )$. The set $A=\beta
(x,\rho _{\mathcal{E}})\cap \beta (y,(4q+1)\lambda )$ is connected and
contained in the union of the tiles of $\mathrm{I}_{\Delta }(\mathcal{E})$.
Consequently, $\mathcal{A}=\{X\in \mathrm{I}_{\Delta }(\mathcal{E})\mid
X\cap A\neq \emptyset \}$\ is connected. Moreover, $T$\ and $W$\ belong to $%
\mathcal{A}$ since $y$ and $z$ belong to $A$.

We are going to prove that, for each $X\in \mathcal{A}$, we have $\mathcal{B}%
_{q}^{\mathcal{E}}(X)\subset \mathcal{E}^{\prime }$ and $Z\varphi =Z\varphi
^{\prime }$\ for each $Z\in \mathcal{B}_{q}^{\mathcal{E}}(X)$.\ As this
property is true for $X=T$, it suffices to show that, if it is true for some 
$X\in \mathcal{A}$, then it is true for any $Y\in \mathcal{A}$ such that $%
X\cap Y\neq \emptyset $.

We have $Y\in \mathcal{B}_{1}^{\mathcal{E}}(X)$, and therefore $\mathcal{B}%
_{q-1}^{\mathcal{E}}(Y)\subset \mathcal{E}^{\prime }$ and $Z\varphi
=Z\varphi ^{\prime }$\ for each $Z\in \mathcal{B}_{q-1}^{\mathcal{E}}(Y)$.\
Moreover, by Lemma 2.7.1, we have $Y\in \mathrm{I}_{\Delta }(\mathcal{E}%
^{\prime })$ since $Y\cap \beta (y,(4q+1)\lambda )\neq \emptyset $, and
therefore $\mathcal{B}_{q}^{\mathcal{E}^{\prime }}(Y)\varphi ^{\prime
}=B_{q}^{M}(Y\varphi ^{\prime })=B_{q}^{M}(Y\varphi )$ since $\varphi
^{\prime }$ is a $\mathcal{L}_{\Delta }$-isomorphism from $\mathcal{E}%
^{\prime }$\ to $B_{r}^{M}(T\varphi )$. We also have $\mathcal{B}_{q}^{%
\mathcal{E}}(Y)\varphi \subset B_{q}^{M}(Y\varphi )$ since $Y\in \mathrm{I}%
_{\Delta }(\mathcal{E})$. Consequently, $\varphi \varphi ^{\prime -1}$ is
defined on $\mathcal{B}_{q}^{\mathcal{E}}(Y)$ and stabilizes the elements of 
$\mathcal{B}_{q-1}^{\mathcal{E}}(Y)$.$\ $It follows $\mathcal{B}_{q}^{%
\mathcal{E}}(Y)=\mathcal{B}_{q}^{\mathcal{E}^{\prime }}(Y)$ and $Z\varphi
=Z\varphi ^{\prime }$ for each $Z\in \mathcal{B}_{q}^{\mathcal{E}}(Y)$.$%
~~\blacksquare $\bigskip

For each relational language $\mathcal{L}$, each $\mathcal{L}$-structure $M$
and each subgroup $H$ of $\mathrm{Aut}(M)$, we define the $\mathcal{L}$%
-structure $M/H$ as follows: The elements of $M/H$ are the classes $xH$ for $%
x\in M$. For $R(u_{1},...,u_{k})\in \mathcal{L}$ and $x_{1},...,x_{k}\in M/H$%
, we write $R(x_{1},...,x_{k})$ if there exist some representatives $%
y_{1},...,y_{k}$ of $x_{1},...,x_{k}$\ in $M$ such that $R(y_{1},...,y_{k})$.

For each $i\in \{1,...,k\}$ and\ each representative $y_{i}$ of $x_{i}$\ in $%
M$, there exist some representatives $y_{1},...,y_{i-1},y_{i+1},...,y_{k}$
of $x_{1},...,x_{i-1},x_{i+1},...,x_{k}$\ in $M$ such that $%
R(y_{1},...,y_{k})$. The canonical surjection from $M$ to $M/H$ is a
homomorphism. If $H$ is normal in $\mathrm{Aut}(M)$, then any automorphism
of $M$\ induces an automorphism of $M/H$.

In Proposition 2.8 below, we do not use the supplementary hypotheses on $%
(E,\delta )$ which were introduced for Theorem 2.7. For each $\Delta $%
\textbf{-}tiling $\mathcal{T}$\ and each subgroup $H$\ of $G$ such that\ $%
T\sigma \in \mathcal{T}-\mathcal{B}_{2q}^{\mathcal{T}}(T)$ for each\ $T\in 
\mathcal{T}$ and each $\sigma \in H$, we denote by $\mathcal{T}/H$\ the
tiling of $E/H$\ induced by $\mathcal{T}$.\ The canonical surjection $\pi :%
\mathcal{T}\rightarrow \mathcal{T}/H$\ is a $\mathcal{L}_{\Delta }$%
-homomorphism. Proposition 2.8 implies that, in Theorem 2.7, the $\mathcal{L}%
_{\Delta }$-structure $M$\ is isomorphic to a quotient of a $\Delta $\textbf{%
-}tiling.\bigskip

\noindent \textbf{Proposition 2.8.} Consider a $\Delta $\textbf{-}tiling $%
\mathcal{T}$, a $\mathcal{L}_{\Delta }$-structure $M$ and a surjective $%
\mathcal{L}_{\Delta }$-homomorphism $\varphi :\mathcal{T}\rightarrow M$
which induces a $\mathcal{L}_{\Delta }$-isomorphism from $\mathcal{B}_{q}^{%
\mathcal{T}}(T)$ to $B_{q}^{M}(T\varphi )$ for each\ $T\in \mathcal{T}$.
Then $\varphi $ induces a $\mathcal{L}_{\Delta }$-isomorphism from $\mathcal{%
T}/H$\ to $M$, where $H=\{\sigma \in G\mid \mathcal{T}\sigma =\mathcal{T}$
and $T\sigma \varphi =T\varphi $ for each\ $T\in \mathcal{T}\}$.\bigskip

\noindent \textbf{Proof}. First we show that $\varphi $ induces a bijective
homomorphism from $\mathcal{T}/H$\ to $M$. For any $S,T\in \mathcal{T}$\
such that $S\varphi =T\varphi $, as $\varphi $\ induces some $\mathcal{L}%
_{\Delta }$-isomorphisms from $\mathcal{B}_{q}^{\mathcal{T}}(S)$ and $%
\mathcal{B}_{q}^{\mathcal{T}}(T)$ to $B_{q}^{M}(S\varphi
)=B_{q}^{M}(T\varphi )$, Theorem 1.3 implies that there exists a unique $%
\sigma _{S,T}\in G$ such that $S\sigma _{S,T}=T$ and $U\varphi =U\sigma
_{S,T}\varphi $ for each $U\in \mathcal{B}_{q}^{\mathcal{T}}(S)$. For any $%
S,T\in \mathcal{T}$\ such that $S\varphi =T\varphi $ and for each $U\in 
\mathcal{B}_{1}^{\mathcal{T}}(S)$, we have $\sigma _{S,T}=\sigma _{U,V}$
where $V=U\sigma _{S,T}$, since $\sigma _{S,T}$ and $\sigma _{U,V}$ coincide
on $\mathcal{B}_{q-1}^{\mathcal{T}}(S)$. Consequently, for any $S,T\in 
\mathcal{T}$\ such that $S\varphi =T\varphi $, we have $\mathcal{T}\sigma
_{S,T}=\mathcal{T}$ and $U\varphi =U\sigma _{S,T}\varphi $ for each $U\in 
\mathcal{T}$.

It remains to be proved that, for each $R(w_{1},...,w_{k})\in \mathcal{L}%
_{\Delta }$ and any $T_{1},...,T_{k}\in \mathcal{T}$ such that $M$ satisfies 
$R(T_{1}\varphi ,...,T_{k}\varphi )$, there exist $U_{1}\in
T_{1}H,...,U_{k}\in T_{k}H$\ such that $\mathcal{T}$\ satisfies $%
R(U_{1},...,U_{k})$. As $\varphi $ induces a $\mathcal{L}_{\Delta }$%
-isomorphism from $\mathcal{B}_{q}^{\mathcal{T}}(T_{1})$ to $%
B_{q}^{M}(T_{1}\varphi )$, there exist $U_{2},...,U_{k}\in \mathcal{B}_{q}^{%
\mathcal{T}}(T_{1})$ such that $R(T_{1},U_{2},...,U_{k})$\ and $U_{2}\varphi
=T_{2}\varphi ,...,U_{k}\varphi =T_{k}\varphi $. We have $U_{2}\in
T_{2}H,...,U_{k}\in T_{k}H$\ according to the first part of the proof.~~$%
\blacksquare $\bigskip

\textbf{3. Periodicity, invariance through a translation.}\bigskip

In the present section, we consider a finite relational language $\mathcal{L}
$. We generalize to uniformly locally finite $\mathcal{L}$-structures the
notions of periodicity and invariance through a nontrivial translation,
which are usually considered for tilings of the euclidean spaces $%
\mathbb{R}
^{n}$. In particular, we obtain generalizations for tilings of noneuclidean
spaces. The notions of mathematical logic used for Proposition 3.1 and
Corollary 3.2 are defined, for instance, in [5].\bigskip

\noindent \textbf{Proposition 3.1.} Consider a formula $\theta (u,v)$ in $%
\mathcal{L}$ and two elementarily equivalent $\mathcal{L}$-structures $M,N$
with $M$ connected locally finite. Suppose that there exist an element $x\in
M$ such that $\left\{ y\in M\mid \theta (x,y)\right\} $ is finite, and an
automorphism $g$ of $N$ such that $\theta (y,yg)$\ for each $y\in N$. Then
there exists an automorphism $f$ of $M$ such that $\theta (y,yf)$\ for each $%
y\in M$. Moreover, for each $r\in 
\mathbb{N}
$, if $\left\{ y\in N\mid yg=y\right\} $ contains no ball $B_{N}(z,r)$, then
we can choose $f$ in such a way that $\left\{ y\in M\mid yf=y\right\} $
contains no ball $B_{M}(z,r)$.\bigskip

\noindent \textbf{Proof}.\ Fix $x\in M$ such that $\left\{ y\in M\mid \theta
(x,y)\right\} $ is finite. For each $k\in 
\mathbb{N}
$, as $B_{M}(x,k)$ is finite, the following property of a $\mathcal{L}$%
-structure $P$ can be expressed by one sentence:

\noindent For each $y\in P\ $such that $(B_{P}(y,k),y)\cong (B_{M}(x,k),x)$,
there exist an element $z\in P$ and an isomorphism $f:(B_{P}(y,k),y)%
\rightarrow (B_{P}(z,k),z)$ such that $\theta (u,uf)$ for each $u\in
B_{P}(y,k)$ (respectively $\theta (u,uf)$ for each $u\in B_{P}(y,k)$ and $%
\left\{ v\in B_{P}(y,k)\mid vf=v\right\} $ contains no ball $B_{P}(u,r)$ for 
$u\in B_{P}(y,k-r)$).

Suppose that there exists an automorphism $g$ of $N$ such that $\theta
(y,yg) $\ for each $y\in N$ (respectively $\theta (y,yg)$\ for each $y\in N$
and $\left\{ y\in N\mid yg=y\right\} $ contains no ball $B_{N}(z,r)$). Then
the sentence considered above is true in $N$, and therefore true in $M$.

For each $k\in 
\mathbb{N}
$, consider the nonempty set $A_{k}$ consisting of the pairs $(y,f)$ with $%
y\in M$ and $f:(B_{M}(x,k),x)\rightarrow (B_{M}(y,k),y)$ isomorphism such
that $\theta (u,uf)$ for each $u\in B_{M}(x,k)$ (respectively $\theta (u,uf)$
for each $u\in B_{M}(x,k)$ and $\left\{ v\in B_{M}(x,k)\mid vf=v\right\} $
contains no ball $B_{M}(u,r)$ for $u\in B_{M}(x,k-r)$). Then each $A_{k}$ is
finite since $M$ is locally finite and $\left\{ y\in M\mid \theta
(x,y)\right\} $ is finite. Moreover, for $0\leq k\leq l$, each pair $%
(z,g)\in A_{l}$ gives by restriction a pair $(y,f)\in A_{k}$.

Consequently, by K\"{o}nig's lemma, there exists a sequence $%
(y_{k},f_{k})_{k\in 
\mathbb{N}
}\in \Pi _{k\in 
\mathbb{N}
}A_{k}$\ with $f_{k}$\ restriction of $f_{l}$\ for $0\leq k\leq l$. The
inductive limit of such a sequence gives an automorphism of $M$ which
satisfies the required properties.~~$\blacksquare $\bigskip

By [7, Theorem 2.3], two locally finite $\mathcal{L}$-structures are
elementarily equivalent if and only if they are locally isomorphic.
Moreover, for any $r,s\in 
\mathbb{N}
^{\ast }$, there exists a formula\ $\theta _{r,s}(u,v)$\ which expresses the
property $v\in B_{N}(u,r)$ in each $\mathcal{L}$-structure $N$ such that $%
\left\vert B_{N}(x,r)\right\vert \leq s$\ for each $x\in N$. Consequently,
we have:\bigskip

\noindent \textbf{Corollary 3.2.} Consider two locally isomorphic\ uniformly
locally finite $\mathcal{L}$-structures $M,N$ which $M$ connected. Let $r\in 
\mathbb{N}
^{\ast }$\ and $s\in 
\mathbb{N}
$.\ Suppose that there exists an automorphism $g$\ of $N$\ such that $yg\in
B_{N}(y,r)$\ for each $y\in N$, and such that $\left\{ y\in N\mid
yg=y\right\} $ contains no ball $B_{N}(z,s)$. Then, there exists an
automorphism $f$ of $M$ such that $yf\in B_{M}(y,r)$\ for each $y\in M$, and
such that $\left\{ y\in M\mid yf=y\right\} $ contains no ball $B_{M}(z,s)$%
.\bigskip

\noindent \textbf{Remark.} In particular, for each $r\in 
\mathbb{N}
^{\ast }$, if $N$ has an automorphism $g$ without fixed point such that $%
yg\in B_{N}(y,r)$ for each $y\in N$, then $M$ has an automorphism $f$
without fixed point such that $yf\in B_{M}(y,r)$ for each $y\in M$.\bigskip

Now we consider the metric space $(E,\delta )$,\ the group $G$\ of
isometries of $E$\ and the set $\Delta $\ defined in Section 1. For each $%
\Delta $-tiling $\mathcal{T}$ and each $\sigma \in G$ such that $\mathcal{T}%
\sigma =\mathcal{T}$, we say that $\sigma $ is a \emph{translation of }$%
\mathcal{T}$\ if there exists $r\in 
\mathbb{N}
^{\ast }$ such that $T\sigma \in \mathcal{B}_{r}^{\mathcal{T}}(T)$\ for each 
$T\in \mathcal{T}$. The set $\mathrm{Trans}(\mathcal{T})$ of all
translations of\emph{\ }$\mathcal{T}$ is a subgroup of $G$.

Any $\sigma \in G$ such that $\mathcal{T}\sigma =\mathcal{T}$ is a
translation of $\mathcal{T}$ if and only if there exists $s\in 
\mathbb{N}
^{\ast }$ such that $T\sigma \in B_{\mathcal{T}}(T,s)$\ for each $T\in 
\mathcal{T}$. Moreover,\ for each $T\in \mathcal{T}$, there exists no $%
\sigma \in G-\{\mathrm{Id}\}$\ such that $S\sigma =S$ for each $S\in B_{%
\mathcal{T}}(T,1)=\mathcal{B}_{2q}^{\mathcal{T}}(T)$. Consequently, the
result below follows from Corollary 3.2:\bigskip

\noindent \textbf{Corollary 3.3}. For any locally isomorphic $\Delta $%
-tilings $\mathcal{S},\mathcal{T}$, we have $\mathrm{Trans}(\mathcal{S}%
)=\left\{ \mathrm{Id}\right\} $\ if and only if $\mathrm{Trans}(\mathcal{T}%
)=\left\{ \mathrm{Id}\right\} $.\bigskip

For each $\Delta $-tiling $\mathcal{T}$ and each $\sigma \in \mathrm{Trans}(%
\mathcal{T})$, $\sup_{x\in E}\delta (x,x\sigma )$ is finite. Conversely, if $%
(E,\delta )$\ is weakly homogeneous, then the second part of Proposition 1.2
implies that any $\sigma \in G$\ such that $\mathcal{T}\sigma =\mathcal{T}$
belongs to $\mathrm{Trans}(\mathcal{T})$ if $\sup_{x\in E}\delta (x,x\sigma
) $ is finite. For each $k\in 
\mathbb{N}
^{\ast }$, if $E$ is the set $%
\mathbb{R}
^{k}$ equipped with a distance defined from a norm, then we have $\sigma \in 
\mathrm{Trans}(\mathcal{T})$ if and only if $\sigma $ is a translation in
the usual sense, since any surjective isometry of $E$ is affine in that case
by [2, Th. 14.1].

The following result generalizes well known properties of the translations
in the euclidean spaces $%
\mathbb{R}
^{k}$. Here, we use the conditions (CVX$\nu $) and (FIX$\nu $) introduced at
the end of Section 1.\bigskip

\noindent \textbf{Theorem 3.4.} Let $\mathcal{T}$ be a $\Delta $-tiling.
Then, for each $\sigma \in \mathrm{Trans}(\mathcal{T})-\left\{ \mathrm{Id}%
\right\} $, there exist $x,y\in E$\ such that $\delta (x,x\sigma
)=\inf_{z\in E}\delta (z,z\sigma )$\ and $\delta (y,y\sigma )=\sup_{z\in
E}\delta (z,z\sigma )$. Moreover, if $(E,\delta )$\ is geodesic and
satisfies (CVX$\nu $) and (FIX$\nu $) for each $\nu \in 
\mathbb{R}
_{>0}$, then $\mathrm{Trans}(\mathcal{T})$ is torsion-free abelian and each $%
\sigma \in \mathrm{Trans}(\mathcal{T})-\left\{ \mathrm{Id}\right\} $ has no
fixed point.\bigskip

\noindent \textbf{Proof.} We fix $\sigma \in \mathrm{Trans}(\mathcal{T}%
)-\left\{ \mathrm{Id}\right\} $, we write $\alpha =\inf_{z\in E}\delta
(z,z\sigma )$ and we show that there exists $x\in E$ such that $\delta
(x,x\sigma )=\alpha $. It can be proved in a similar way that there exists $%
y\in E$ such that $\delta (y,y\sigma )=\sup_{z\in E}\delta (z,z\sigma )$.

We consider $r\in 
\mathbb{N}
^{\ast }$ such that $T\sigma \in \mathcal{B}_{r}^{\mathcal{T}}(T)$ for each $%
T\in \mathcal{T}$. For each $T\in \mathcal{T}$, $\sigma $ induces a
bijection $\sigma _{T}:\mathcal{B}_{q}^{\mathcal{T}}(T)\rightarrow \mathcal{B%
}_{q}^{\mathcal{T}}(T\sigma )$.\ The triples $(\mathcal{B}_{q+r}^{\mathcal{T}%
}(T),T,\sigma _{T})$\ for $T\in \mathcal{T}$ fall in finitely many
isomorphism classes. Consequently, there exist two sequences $(x_{k})_{k\in 
\mathbb{N}
}\subset E$ and $(T_{k})_{k\in 
\mathbb{N}
}\subset \mathcal{T}$\ with $x_{k}\in T_{k}$\ for each $k\in 
\mathbb{N}
$\ such that $\lim \delta (x_{k},x_{k}\sigma )=\alpha $ and such that all
the triples $(\mathcal{B}_{q+r}^{\mathcal{T}}(T_{k}),T_{k},\sigma _{T_{k}})$%
\ are isomorphic.

For each $k\in 
\mathbb{N}
$, we consider $\tau _{k}\in G$\ which induces an isomorphism from $(%
\mathcal{B}_{q+r}^{\mathcal{T}}(T_{0}),T_{0},\sigma _{T_{0}})$ to $(\mathcal{%
B}_{q+r}^{\mathcal{T}}(T_{k}),T_{k},\sigma _{T_{k}})$. We have $\sigma =\tau
_{k}\sigma \tau _{k}^{-1}$ for each $k\in 
\mathbb{N}
$\ since $T\sigma =T\sigma _{T_{0}}=T\tau _{k}\sigma _{T_{k}}\tau
_{k}^{-1}=T\tau _{k}\sigma \tau _{k}^{-1}$\ for each $T\in \mathcal{B}_{q}^{%
\mathcal{T}}(T_{0})$. Consequently, the elements $y_{k}=x_{k}\tau
_{k}^{-1}\in T_{0}$\ satisfy $\delta (y_{k},y_{k}\sigma )=\delta
(x_{k},x_{k}\sigma )$, and therefore $\lim \delta (y_{k},y_{k}\sigma
)=\alpha $. As $T_{0}$\ is compact, it follows that there exists $x\in T_{0}$%
\ such that $\delta (x,x\sigma )=\alpha $.

Now suppose that $(E,\delta )$\ is geodesic and satisfies (CVX$\nu $) and
(FIX$\nu $) for each $\nu \in 
\mathbb{R}
_{>0}$. Consider $\sigma \in \mathrm{Trans}(\mathcal{T})$ and $s\in 
\mathbb{N}
-\{0,1\}$ such that $\sigma ^{s}=\mathrm{Id}$. Then, for each $x\in E$, we
have $A_{x}\sigma =A_{x}$\ for $A_{x}=\{x,x\sigma ,...,x\sigma ^{s-1}\}$; by
Lemma 1.9, there exists a unique $w_{x}\in E$\ such that $\delta
(w_{x},z)\leq \mathrm{Rad}(A_{x})$ for each $z\in A_{x}$; it follows $%
w_{x}\sigma =w_{x}$. The properties (FIX$\nu $) imply $\sigma =\mathrm{Id}$
since $\delta (x,w_{x})\leq \mathrm{Rad}(A_{x})\leq (s/2)\sup_{z\in E}\delta
(z,z\sigma )$ for each $x\in E$.

Now consider any $\sigma \in \mathrm{Trans}(\mathcal{T})$ with a fixed point 
$x$. Then, for each $y\in E$, there exists $s\in 
\mathbb{N}
^{\ast }$ such that $y\sigma ^{s}=y$, since $\mathrm{Trans}(\mathcal{T})$ is
discrete by Proposition 1.8 and the elements $y\sigma ^{k}$ for $k\in 
\mathbb{Z}
$ all belong to $\beta (x,\delta (x,y))$. Consider $\zeta ,\eta \in 
\mathbb{R}
_{>0}$ such that no $\tau \in G-\{\mathrm{Id}\}$ fixes the points of a set $%
A\subset E$\ with $\beta (y,\eta )\subset \cup _{z\in A}\beta (z,\zeta )$
for some $y\in E$, and choose a finite set $A$ with that property. Then
there exists $s\in 
\mathbb{N}
^{\ast }$ such that $z\sigma ^{s}=z$ for each $z\in A$, which implies $%
\sigma ^{s}=\mathrm{Id}$ and $\sigma =\mathrm{Id}$.

It remains to be proved that any $\sigma ,\tau \in \mathrm{Trans}(\mathcal{T}%
)$ commute.\ It suffices to show that $\sigma $\ commutes with $\tau ^{r}$
for some $r\in 
\mathbb{Z}
^{\ast }$, since $(\sigma ^{-1}\tau \sigma )^{r}=\sigma ^{-1}\tau ^{r}\sigma
=\mathrm{Id}$ implies $\sigma ^{-1}\tau \sigma =\mathrm{Id}$.

We fix $x\in E$\ and we consider $\alpha \in 
\mathbb{R}
_{>0}$ such that $\delta (y,y\sigma )\leq \alpha $ for each $y\in E$. For
each $r\in 
\mathbb{Z}
$, we have $\delta (x,x\sigma ^{-1}\tau ^{-r}\sigma \tau ^{r})\leq 2\alpha $%
\ since $\delta (x,x\sigma ^{-1})\leq \alpha $\ and $\delta (x\sigma
^{-1},x\sigma ^{-1}\tau ^{-r}\sigma \tau ^{r})=\delta (x\sigma ^{-1}\tau
^{-r},x\sigma ^{-1}\tau ^{-r}\sigma )\leq \alpha $.

As $\mathrm{Trans}(\mathcal{T})$ is discrete, it follows that there exist $%
r\neq s$ in $%
\mathbb{Z}
$\ such that $x\sigma ^{-1}\tau ^{-r}\sigma \tau ^{r}=x\sigma ^{-1}\tau
^{-s}\sigma \tau ^{s}$. We have $\sigma ^{-1}\tau ^{-r}\sigma \tau
^{r}=\sigma ^{-1}\tau ^{-s}\sigma \tau ^{s}$\ since the elements of $\mathrm{%
Trans}(\mathcal{T})$ have no fixed point. It follows $\tau ^{-r}\sigma \tau
^{r}=\tau ^{-s}\sigma \tau ^{s}$\ and $\tau ^{s-r}\sigma =\sigma \tau ^{s-r}$%
.~~$\blacksquare $\bigskip

Now we generalize the notion of periodicity to connected $\mathcal{L}$%
-structures. We say that such a structure $M$ is \emph{periodic} if it
contains a finite set $A$\ such that $M=\cup _{\sigma \in \mathrm{Aut}%
(M)}A\sigma $.\bigskip

\noindent \textbf{Proposition 3.5.} Let $M,N$ be locally isomorphic locally
finite connected $\mathcal{L}$-structures. If $N$ is periodic then $M$ is
isomorphic to $N$.\bigskip

\noindent \textbf{Proof.} We fix $w\in M$. As $N$\ is periodic and
connected, there exist $z\in N$ and $r\in 
\mathbb{N}
$\ such that $N=\cup _{\sigma \in \mathrm{Aut}(N)}B_{N}(z,r)\sigma $. For
each $s\in 
\mathbb{N}
$, as $M$ and $N$ are locally isomorphic, there exists $x_{s}\in N$\ such
that $(B_{M}(w,s),w)\cong (B_{N}(x_{s},s),x_{s})$, and therefore $y_{s}\in
B_{N}(z,r)$\ such that $(B_{M}(w,s),w)\cong (B_{N}(y_{s},s),y_{s})$. We
consider the nonempty set $A_{s}$ consisting of the isomorphisms $\theta
:(B_{M}(w,s),w)\rightarrow (B_{N}(y,s),y)$ with $y\in B_{N}(z,r)$.

The sets $A_{s}$ are finite since $N$ is locally finite. Moreover, for any $%
s\leq t$, the restriction of each $\theta \in A_{t}$\ to $(B_{M}(w,s),w)$
belongs to $A_{s}$.\ Consequently, by K\"{o}nig's Lemma, there exists a
strictly increasing sequence $(\theta _{s})_{s\in 
\mathbb{N}
}$ with $\theta _{s}\in A_{s}$ for each $s\in 
\mathbb{N}
$. The inductive limit is an isomorphism $\theta :(M,w)\rightarrow (N,y)$\
with $y\in B_{N}(z,r)$.~~$\blacksquare $\bigskip

We call a \emph{period} of a $\mathcal{L}$-structure $M$ any set $A\subset M$%
\ such that $M$ is the disjoint union of the sets $A\sigma $\ for $\sigma
\in \mathrm{Aut}(M)$. If $M$ is periodic and if it has some periods, then
all of them have the same finite number of elements. We call it the \emph{%
periodicity rank} of $M$.

For each $\mathcal{L}$-structure $M$, we say that $A\subset M$\ is \emph{%
weakly connected} if, for each subset $B$ with $\varnothing \subsetneq
B\subsetneq A$, there exist $k\in 
\mathbb{N}
-\{0,1\}$\ and $R(u_{1},...,u_{k})\in \mathcal{L}$\ which is satisfied by
some $(x_{1},...,x_{k})\in M^{k}$\ with $\{x_{1},...,x_{k}\}\cap B\neq
\varnothing $\ and $\{x_{1},...,x_{k}\}\cap (A-B)\neq \varnothing $.\bigskip

\noindent \textbf{Proposition 3.6.} Any connected $\mathcal{L}$-structure $M$
has a weakly connected period if the nontrivial automorphisms of $M$ have no
fixed point.\bigskip

\noindent \textbf{Proof.} We show that $M=\cup _{\sigma \in \mathrm{Aut}%
(M)}A\sigma $ for each $A\subset M$ which is maximal for the conjunction of
the two properties: $A$ weakly connected\ and $A\cap A\sigma =\varnothing $\
for each $\sigma \in \mathrm{Aut}(M)-\{\mathrm{Id}\}$.

Otherwise, as $M$ is connected, there exist $k\in 
\mathbb{N}
-\{0,1\}$, $(x_{1},...,x_{k})\in M^{k}$ which satisfies some $%
R(u_{1},...,u_{k})\in \mathcal{L}$ and $1\leq i,j\leq k$ such that $x_{i}\in
\cup _{\sigma \in \mathrm{Aut}(M)}A\sigma $\ and $x_{j}\in M-\cup _{\sigma
\in \mathrm{Aut}(M)}A\sigma $. We consider $\tau \in \mathrm{Aut}(M)$\ such
that $x_{i}\tau \in A$. Then $B=A$ $\cup $ $\left\{ x_{j}\tau \right\} $\ is
weakly connected since $(x_{1}\tau ,...,x_{k}\tau )$ satisfies $R$. As $%
x_{j}\tau $ is not a fixed point of a nontrivial automorphism of $M$ and
does not belong to $\cup _{\sigma \in \mathrm{Aut}(M)}A\sigma $, we have $%
B\sigma \cap B=\varnothing $\ for each $\sigma \in \mathrm{Aut}(M)-\{\mathrm{%
Id}\}$, contrary to the maximality of $A$.~~$\blacksquare $\bigskip

Now we introduce some supplementary conditions which will be used for the
investigation of periodic $\mathcal{L}$-structures.

We say that a $\mathcal{L}$-structure $M$ is \emph{equational} (cf. [7,
Section 4]) if $R(x_{1},...,x_{k})$, $R(y_{1},...,y_{k})$ and $x_{i}=y_{i}$
imply $x_{j}=y_{j}$ for each $k\in 
\mathbb{N}
^{\ast }$, each $R(u_{1},...,u_{k})\in \mathcal{L}$, any $%
(x_{1},...,x_{k}),(y_{1},...,y_{k})\in M^{k}$\ and any $i,j\in \{1,...,k\}$.

Any equational $\mathcal{L}$-structure is uniformly locally finite. If it is
connected, then its nontrivial automorphisms have no fixed point.

For each equational $\mathcal{L}$-structure $M$, each $k\in 
\mathbb{N}
-\{0,1\}$, each $R(u_{1},...,u_{k})\in \mathcal{L}$, any $i,j\in \{1,...,k\}$%
\ such that $i\neq j$\ and any $x,y\in M$, we write $x(R,i,j)=y$ if there
exists $(z_{1},...,z_{k})\in M^{k}$ which satisfies $R$\ with $z_{i}=x$\ and 
$z_{j}=y$.

For each $n\in 
\mathbb{N}
$, each word $w=(R_{1},i_{1},j_{1})...(R_{n},i_{n},j_{n})$\ and any $x,y\in
M $, we write $xw=y$\ if there exist $z_{0},...,z_{n}\in M$\ such that $%
z_{0}=x $, $z_{n}=y$\ and $z_{m-1}(R_{m},i_{m},j_{m})=z_{m}$\ for $1\leq
m\leq n$.\ We denote by $\Omega _{\mathcal{L}}$\ the set of all such words.

For any equational $\mathcal{L}$-structures $M,N$, each $x\in M$ and each $%
r\in 
\mathbb{N}
^{\ast }$, we call a \emph{partial isomorphism} from $B_{M}(x,r)$ to $N$ any
injective map $\rho $\ such that, for each $k\in 
\mathbb{N}
^{\ast }$, each $R(u_{1},...,u_{k})\in \mathcal{L}$, any $i,j\in \{1,...,k\}$%
\ and any $y,z\in B_{M}(x,r)$, $y(R,i,j)$ exists if and only if $y\rho
(R,i,j)$ exists and $y(R,i,j)=z$ if and only if $y\rho (R,i,j)=z\rho $.

If $\rho $\ is a partial isomorphism from $B_{M}(x,r)$ to $N$, then $\rho $\
is an isomorphism from $(B_{M}(x,r),x)$ to $(B_{N}(x\rho ,r),x\rho )$ and $%
\rho ^{-1}$\ is a partial isomorphism from\ $B_{N}(x\rho ,r)$\ to $M$. Any
isomorphism $\sigma :(B_{M}(x,r+1),x)\rightarrow (B_{N}(x\sigma
,r+1),x\sigma )$ gives by restriction a partial isomorphism from $B_{M}(x,r)$
to $N$.

We say that an equational $\mathcal{L}$-structure $M$ is \emph{commutative}
if, for each $x\in M$ and each $r\in 
\mathbb{N}
^{\ast }$, each partial isomorphism $\rho :B_{M}(x,r)\rightarrow M$
satisfies $y\rho \sigma =y\sigma \rho $ for each automorphism $\sigma $ of $%
M $ and each $y\in B_{M}(x,r)$\ such that $y\sigma \in B_{M}(x,r)$.

We say that $M$ is \emph{strongly commutative} if we have $xvw=xwv$ for each 
$x\in M$\ and any $v,w\in \Omega _{\mathcal{L}}$\ such that $xvw$ and $xwv$
exist. Any strongly commutative connected equational $\mathcal{L}$-structure
is commutative since, for $x,y,\rho ,\sigma $ defined as above, there exist $%
v,w\in \Omega _{\mathcal{L}}$\ such that the equalities $y\rho =yv$ and $%
y\sigma =yw$ are respectively true in $M$\ and in $B_{M}(x,r)$, which
implies $y\rho \sigma =yvw$ and $y\sigma \rho =ywv$.

We say that a class $\mathcal{E}$\ of equational $\mathcal{L}$-structures is 
\emph{strongly regular} (cf. [7, Section 4]) if $xw=x$ and $yw=y$ are
equivalent for any $M,N\in \mathcal{E}$, each $x\in M$, each $y\in N$\ and
each $w\in \Omega _{\mathcal{L}}$ such that $xw$ and $yw$ exist.

Equationality, strong commutativity and strong regularity are local
properties. They are preserved by local isomorphism.\bigskip

\noindent \textbf{Proposition 3.7.} Let $M$ be a connected equational
commutative $\mathcal{L}$-structure. If $M$ is periodic of rank $r\in 
\mathbb{N}
^{\ast }$, then, for each $x\in M$ and each integer $s\geq r$, each partial
isomorphism from $B_{M}(x,s)$ to $M$ can be extended into a unique
automorphism of $M$.\bigskip

\noindent \textbf{Proof.} The extension is necessarily unique since the
nontrivial automorphisms of an equational $\mathcal{L}$-structure have no
fixed point. It remains to be proved that it exists.

By Proposition 3.6, $M$ has a weakly connected period $A$ with $\left\vert
A\right\vert =r$. Replacing $A$ by its image through an appropriate
automorphism of $M$, we reduce the proof to the case $x\in A$. Then $A$ is
contained in $B_{M}(x,r-1)$.

For each partial isomorphism $\rho :B_{M}(x,s)\rightarrow M$, we define a
map $\overline{\rho }:M\rightarrow M$\ as follows:\ for each $y\in M$, as
the nontrivial automorphisms of $M$ have no fixed point, there exists a
unique $\sigma \in \mathrm{Aut}(M)$\ such that $y\sigma \in A$; we write $y%
\overline{\rho }=y\sigma \rho \sigma ^{-1}$. For each $y\in B_{M}(x,s)$, we
have $y\overline{\rho }=y\rho $ because the commutativity of $M$ implies $%
y\rho \sigma =y\sigma \rho $. We observe that $\overline{\rho }^{-1}$ is
defined in the same way from the partial isomorphism $\rho ^{-1}:B_{M}(x\rho
,s)\rightarrow M$\ and the weakly connected period $A\rho $. Consequently, $%
\overline{\rho }$ is bijective and, by reason of symmetry, it suffices to
show that $\overline{\rho }$\ is a homomorphism.

For each $k\in 
\mathbb{N}
^{\ast }$, each $R(u_{1},...,u_{k})\in \mathcal{L}$ and each $%
(x_{1},...,x_{k})\in M^{k}$\ which satisfies $R$,\ we consider $\sigma
_{1},...,\sigma _{k}\in \mathrm{Aut}(M)$\ such that $x_{1}\sigma
_{1},...,x_{k}\sigma _{k}$ belong to $A$. Then $x_{1}\sigma
_{1},...,x_{k}\sigma _{k}$ belong to $B_{M}(x,r-1)$ and $x_{2}\sigma
_{1},...,x_{k}\sigma _{1}$ belong to $B_{M}(x,r)$ since $(x_{1}\sigma
_{1},x_{2}\sigma _{1},...,x_{k}\sigma _{1})$ satisfies $R$.

We have $x_{1}\overline{\rho }=x_{1}\sigma _{1}\rho \sigma _{1}^{-1}$. For $%
2\leq i\leq k$, we have $x_{i}\overline{\rho }=x_{i}\sigma _{i}\rho \sigma
_{i}^{-1}=(x_{i}\sigma _{1})(\sigma _{1}^{-1}\sigma _{i})\rho (\sigma
_{1}^{-1}\sigma _{i})^{-1}\sigma _{1}^{-1}$. As $x_{i}\sigma _{1}$\ and $%
(x_{i}\sigma _{1})(\sigma _{1}^{-1}\sigma _{i})=x_{i}\sigma _{i}$ belong to $%
B_{M}(x,r)$, the commutativity of $M$ implies

\noindent $(x_{i}\sigma _{1})(\sigma _{1}^{-1}\sigma _{i})\rho =(x_{i}\sigma
_{1})\rho (\sigma _{1}^{-1}\sigma _{i})$, and therefore

\noindent $x_{i}\overline{\rho }=(x_{i}\sigma _{1})(\sigma _{1}^{-1}\sigma
_{i})\rho (\sigma _{1}^{-1}\sigma _{i})^{-1}\sigma _{1}^{-1}=(x_{i}\sigma
_{1})\rho (\sigma _{1}^{-1}\sigma _{i})(\sigma _{1}^{-1}\sigma
_{i})^{-1}\sigma _{1}^{-1}=x_{i}\sigma _{1}\rho \sigma _{1}^{-1}$.

Moreover, $x_{1}\sigma _{1}\rho ,x_{2}\sigma _{1}\rho ,...,x_{k}\sigma
_{1}\rho $ are all defined since $x_{1}\sigma _{1},x_{2}\sigma
_{1},...,x_{k}\sigma _{1}$ belong to $B_{M}(x,r)$. Consequently, $%
(x_{1}\sigma _{1}\rho ,x_{2}\sigma _{1}\rho ,...,x_{k}\sigma _{1}\rho )$
satisfies $R$ like $(x_{1}\sigma _{1},x_{2}\sigma _{1},...,x_{k}\sigma _{1})$
and

\noindent $(x_{1}\overline{\rho },x_{2}\overline{\rho },...,x_{k}\overline{%
\rho })=(x_{1}\sigma _{1}\rho \sigma _{1}^{-1},x_{2}\sigma _{1}\rho \sigma
_{1}^{-1},...,x_{k}\sigma _{1}\rho \sigma _{1}^{-1})$ also$\ $satisfies $R$%
.~~$\blacksquare $\bigskip

\noindent \textbf{Lemma 3.8.} Let $M,N$ be equational $\mathcal{L}$%
-structures with $N$ connected, commutative and periodic of rank $r\in 
\mathbb{N}
^{\ast }$. For each $x\in M$ and each integer $s\geq r$, if there exists a
partial isomorphism from $B_{M}(x,s+1)$ to $N$,\ then each partial
isomorphism from $B_{M}(x,s)$ to $N$ can be extended into a partial
isomorphism from $B_{M}(x,s+1)$ to $N$.\bigskip

\noindent \textbf{Proof.} Let us consider a partial isomorphism $\rho
:B_{M}(x,s)\rightarrow N$ and a partial isomorphism $\sigma
:B_{M}(x,s+1)\rightarrow N$. Then $\rho ^{-1}\sigma $\ is a partial
isomorphism from $B_{N}(x\rho ,s)$ to $N$\ which\ can be extended into a
unique automorphism $\theta $ of $N$ by Proposition 3.7, and $\sigma \theta
^{-1}$\ is a partial isomorphism from $B_{M}(x,s+1)$ to $N$\ which extends $%
\rho $.~~$\blacksquare $\bigskip

\noindent \textbf{Lemma 3.9.} Let $M,N$ be equational $\mathcal{L}$%
-structures with $\left\{ M,N\right\} $ strongly regular. Then,\ for each $%
x\in M$ and each $r\in 
\mathbb{N}
^{\ast }$,\ each partial isomorphism $\rho :B_{M}(x,r)\rightarrow N$ can be
extended into a partial isomorphism from $B_{M}(x,r+1)$ to $N$\ if, for each 
$y\in B_{M}(x,r)$, there exists a partial isomorphism $\rho
_{y}:B_{M}(y,1)\rightarrow N$ such that $y\rho _{y}=y\rho $.\bigskip

\noindent \textbf{Proof.} We define a map $\sigma :B_{M}(x,r+1)\rightarrow N$
as follows:\ for each $z\in B_{M}(x,r+1)$, we consider $y\in B_{M}(x,r)$
such that $z\in B_{M}(y,1)$ and we write $z\sigma =z\rho _{y}$.

For any $z_{1},z_{2}\in B_{M}(x,r+1)$ and any $y_{1},y_{2}\in B_{M}(x,r)$\
such that $z_{1}\in B_{M}(y_{1},1)$ and $z_{2}\in B_{M}(y_{2},1)$, there
exist $v,(R_{1},i_{1},j_{1}),(R_{2},i_{2},j_{2})\in \Omega _{\mathcal{L}}$\
such that the equalities $y_{1}v=y_{2}$, $y_{1}(R_{1},i_{1},j_{1})=z_{1}$
and $y_{2}(R_{2},i_{2},j_{2})=z_{2}$\ are respectively true in the connected 
$\mathcal{L}$-structures $B_{M}(x,r)$, $B_{M}(y_{1},1)$ and $B_{M}(y_{2},1)$%
. We have $z_{2}=z_{1}w$\ for $w=(R_{1},j_{1},i_{1})v(R_{2},i_{2},j_{2})$.
We also have $z_{2}\rho _{y_{2}}=z_{1}\rho _{y_{1}}w$\ since the equalities $%
z_{1}\rho _{y_{1}}(R_{1},j_{1},i_{1})=y_{1}\rho _{y_{1}}=y_{1}\rho $, $%
y_{1}\rho v=y_{2}\rho =y_{2}\rho _{y_{2}}$ and $y_{2}\rho
_{y_{2}}(R_{2},i_{2},j_{2})=z_{2}\rho _{y_{2}}$\ are respectively true in $%
B_{N}(y_{1}\rho _{y_{1}},1)$, $B_{N}(x\rho ,r)$ and $B_{N}(y_{2}\rho
_{y_{2}},1)$.

By strong regularity, the equalities $z_{2}=z_{1}w$\ and $z_{2}\rho
_{y_{2}}=z_{1}\rho _{y_{1}}w$\ imply that $z_{1}\rho _{y_{1}}=z_{2}\rho
_{y_{2}}$ if and only if $z_{1}=z_{2}$. Consequently, $\sigma $ is injective
and the definition of $z\sigma $ given above does not depend on the choice
of the element $y\in B_{M}(x,r)$ such that $z\in B_{M}(y,1)$. It follows
that $z\sigma =z\rho $\ for each $z\in B_{M}(x,r)$ and $z\sigma =z\rho _{y}$%
\ for each $y\in B_{M}(x,r)$ and each $z\in B_{M}(y,1)$. For each $z\in
B_{M}(x,r+1)$ and each $(R,i,j)\in \Omega _{\mathcal{L}}$, the element $%
z\sigma (R,i,j)$ exists if and only if $z(R,i,j)$ exists, because, for each $%
y\in B_{M}(x,r)$, $z\rho _{y}(R,i,j)$ exists if and only if $z(R,i,j)$
exists.

Now let us consider again $y_{1},y_{2},z_{1},z_{2},w$ as above. For each$%
(R,i,j)$ such that $z_{1}(R,i,j)$ exists, or equivalently such that $%
z_{1}\sigma (R,i,j)$ exists, we have $z_{1}(R,i,j)w^{\prime }=z_{2}$ and $%
z_{1}\sigma (R,i,j)w^{\prime }=z_{2}\sigma $ for $w^{\prime }=(R,j,i)w$.
Strong regularity implies that $z_{1}(R,i,j)=z_{2}$ if and only if $%
z_{1}\sigma (R,i,j)=z_{2}\sigma $.~~$\blacksquare $\bigskip

\noindent \textbf{Theorem 3.10.} Let $M,N$ be connected periodic $\mathcal{L}
$-structures of ranks $r,s\geq 1$ with $N$ commutative and $\left\{
M,N\right\} $ strongly regular. Then,\ for each $x\in M$,\ each partial
isomorphism from $B_{M}(x,r+s)$ to $N$ can be extended into an isomorphism
from $M$ to $N$.\bigskip

\noindent \textbf{Proof.} Otherwise, there exist an integer $t\geq r+s$\ and
a partial isomorphism $\rho :B_{M}(x,t)\rightarrow N$ which cannot be
extended into a partial isomorphism from $B_{M}(x,t+1)$ to $N$. By lemma
3.9, there exists\ $y\in B_{M}(x,t)-B_{M}(x,t-1)$ such that no partial
isomorphism $\sigma :B_{M}(y,1)\rightarrow N$ satisfies $y\sigma =y\rho $.

For each $k\in \{0,...,r\}$, we consider $z_{k}\in
B_{M}(x,t-s-k)-B_{M}(x,t-s-k-1)$\ such that $y\in B_{M}(z_{k},s+k)$. The
restriction $\rho _{k}$\ of $\rho $\ to $B_{M}(z_{k},s+k)$\ cannot be
extended into a partial isomorphism from $B_{M}(z_{k},s+k+1)$ to $N$.

As $M$ is periodic of rank $r$, there exist two integers $0\leq
k_{1}<k_{2}\leq r$ and an isomorphism $\sigma :(M,z_{k_{1}})\rightarrow
(M,z_{k_{2}})$. Then $\sigma \rho _{k_{2}}$\ is a partial isomorphism from $%
B_{M}(z_{k_{1}},s+k_{2})$\ to $N$. By Lemma 3.8, it follows that $\rho
_{k_{1}}$\ can be extended into a partial isomorphism from $%
B_{M}(z_{k_{1}},s+k_{2})$ to $N$, whence a contradiction.~~$\blacksquare $%
\bigskip

\noindent \textbf{Theorem 3.11.} For each set $\Sigma $ of local rules for $%
\mathcal{L}$, if the connected $\mathcal{L}$-structures satisfying $\Sigma $
are periodic, equational and commutative, and if their class $\mathcal{E}$
is strongly regular, then $\mathcal{E}$ is a finite union of isomorphism
classes.\bigskip

\noindent \textbf{Proof.} Otherwise, it follows from Theorem 3.10 that there
exists a sequence $(N_{i})_{i\in 
\mathbb{N}
}$\ of periodic $\mathcal{L}$-structures with strictly increasing
periodicity ranks which satisfy $\Sigma $. For each $q\in 
\mathbb{N}
^{\ast }$, the pairs $(B_{N_{i}}(y,q),y)$ for $i\in 
\mathbb{N}
$\ and $y\in N_{i}$\ fall in finitely many isomorphism classes. K\"{o}nig's
lemma applied to the isomorphism classes of pairs $(B_{N_{i}}(y,q),y)$ for $%
q\in 
\mathbb{N}
^{\ast }$, $i\geq q$\ and $y\in N_{i}$\ implies that there exist two
strictly increasing sequences of integers $(k(i))_{i\in 
\mathbb{N}
}$\ and $(q_{i})_{i\in 
\mathbb{N}
}$, and a sequence $(y_{i})_{i\in 
\mathbb{N}
}\in \times _{i\in 
\mathbb{N}
}N_{k(i)}$\ such that $(B_{N_{k(i)}}(y_{i},q_{i}),y_{i})\cong
(B_{N_{k(j)}}(y_{j},q_{i}),y_{j})$\ for $i\leq j$.

Consequently, there exist a sequence $(M_{i})_{i\in 
\mathbb{N}
}$\ of periodic $\mathcal{L}$-structures with strictly increasing
periodicity ranks which satisfy $\Sigma $ and a sequence $(x_{i})_{i\in 
\mathbb{N}
}\in \times _{i\in 
\mathbb{N}
}M_{i}$ such that $(B_{M_{i}}(x_{i},i),x_{i})\cong
(B_{M_{j}}(x_{j},i),x_{j}) $\ for $i\leq j$. The inductive limit of the
pairs $(B_{M_{i}}(x_{i},i),x_{i})$ for these isomorphisms, which are unique
and compatible because of equationality, is a pair $(M,x)$ with $M$
satisfying $\Sigma $ and $x\in M$. We denote by $r$ the periodicity rank of $%
M$ and, for each $i\in 
\mathbb{N}
$, $r_{i}$ the periodicity rank of $M_{i}$.

For each $i\geq r$, we have $r_{i}>r$. Consequently, $M_{i}$ is not
isomorphic to $M$. As $(B_{M_{i}}(x_{i},i),x_{i})$ and $(B_{M}(x,i),x)$\ are
isomorphic, there exist an integer $s_{i}\geq i-1$ and a partial isomorphism 
$\rho _{i}:B_{M_{i}}(x_{i},s_{i})\rightarrow M$\ which cannot be extended
into a partial isomorphism from $B_{M_{i}}(x_{i},s_{i}+1)$ to $M$. By Lemma
3.9, there exists $y_{i}\in B_{M_{i}}(x_{i},s_{i})-B_{M_{i}}(x_{i},s_{i}-1)$%
\ such that no partial isomorphism $\sigma :B_{M_{i}}(y_{i},1)\rightarrow M$
satisfies $y_{i}\sigma =y_{i}\rho _{i}$.

For any integers $i\geq r$ and $1\leq s\leq s_{i}$, we consider $z_{i,s}\in
B_{M_{i}}(x_{i},s_{i}-s)$ such that $y_{i}\in B_{M_{i}}(z_{i,s},s)$. Then\ $%
B_{M_{i}}(y_{i},2s)$\ contains\ $B_{M_{i}}(z_{i,s},s)$.

It follows from K\"{o}nig's Lemma applied to the pairs $%
(B_{M_{i}}(y_{i},2s+1),y_{i})$ for $i\geq r$ and $1\leq s\leq s_{i}$\ that
there exist two strictly increasing sequences of integers $(k(i))_{i\in 
\mathbb{N}
}$\ and $(t_{i})_{i\in 
\mathbb{N}
}$, with $k(i)\geq r$\ and $t_{i}\leq s_{k(i)}$\ for each $i\in 
\mathbb{N}
$, such that $(B_{M_{k(i)}}(y_{k(i)},2t_{i}+1),y_{k(i)})\cong
(B_{M_{k(j)}}(y_{k(j)},2t_{i}+1),y_{k(j)})$\ for $i\leq j$. We denote by $%
(N,y)$ the inductive limit of the pairs $%
(B_{M_{k(i)}}(y_{k(i)},2t_{i}+1),y_{k(i)})$ relative to these isomorphisms,
which are unique and compatible because of equationality. The $\mathcal{L}$%
-structure $N$ is periodic because it satisfies $\Sigma $. We denote by $t$
the periodicity rank of $N$.

For each $i\in 
\mathbb{N}
$, we have $y\in B_{N}(z_{i},t_{i})$ where $z_{i}$\ is the image of $%
z_{k(i),t_{i}}$ in $N$. Moreover $\rho _{k(i)}$ induces a partial
isomorphism $\sigma _{i}:B_{N}(z_{i},t_{i})\rightarrow M$ such that no
partial isomorphism $\tau :B_{N}(y,1)\rightarrow M$ satisfies $y\tau
=y\sigma _{i}$. For $i\geq r+t-1$, this property contradicts Theorem 3.10
since we have $t_{i}\geq r+t$.~~$\blacksquare $\bigskip

Now we consider again the metric space $(E,\delta )$,\ the group $G$\ of
isometries of $E$\ and the set $\Delta $\ defined in Section 1. We say that
a $\Delta $-tiling $\mathcal{T}$ is \emph{periodic} if there exists a finite
subset $\mathcal{E}$ of $\mathcal{T}$ such that $\mathcal{T}$ is the union
of the subsets $\mathcal{E}\sigma $\ for $\sigma \in G$\ such that $\mathcal{%
T}\sigma =\mathcal{T}$. By Theorem 1.3, this property is true if and only if
the $\mathcal{L}_{\Delta }$-structure associated to $\mathcal{T}$ is
periodic.

For each $k\in 
\mathbb{N}
^{\ast }$, if $(E,\delta )$ is the space $%
\mathbb{R}
^{k}$ equipped with a distance defined from a norm and if $G$ consists of
the translations of $E$, then our notion of periodicity coincides with the
classical one. Consequently, the following result generalizes [1, Th. 3.8]
which was proved for tilings of $%
\mathbb{R}
^{2}$ by square tiles:\bigskip

\noindent \textbf{Theorem 3.12.} If $(E,\delta )$ is geodesic and satisfies
(CVX$\lambda $) where $\lambda $\ is the maximum of the radii of the tiles,
if $G$ is commutative and if the elements of $G-\{\mathrm{Id}\}$ have no
fixed point, then each class of periodic $\Delta $-tilings defined by local
rules is a finite union of isomorphism classes.\bigskip

\noindent \textbf{Lemma 3.12.1.} The $\mathcal{L}_{\Delta }$-structures
associated to $\Delta $-tilings are equational, strongly commutative, and
form a strongly regular class.\bigskip

\noindent \textbf{Proof of Lemma 3.12.1.} If some $\sigma \in G$\ stabilizes
a tile, then $\sigma =\mathrm{Id}$ since $\sigma $\ has a fixed point by
Lemma 1.9.\ Consequently, Theorem 1.3 implies that the $\mathcal{L}_{\Delta
} $-structures associated to $\Delta $-tilings are equational. Moreover, for
any $\Delta $-tilings $\mathcal{S},\mathcal{T}$, each $(R,i,j)\in \Omega _{%
\mathcal{L}}$, each $S\in \mathcal{S}$\ and each $T\in \mathcal{T}$, if $%
S(R,i,j)$ and $T(R,i,j)$\ exist, then there exists a unique $\sigma \in G$\
such that $S\sigma =T$ and it satisfies $S(R,i,j)\sigma =T(R,i,j)$.\ By
induction, it follows that, for each $w\in \Omega _{\mathcal{L}}$, each $%
S\in \mathcal{S}$\ and each $T\in \mathcal{T}$, if $Sw$ and $Tw$\ exist,
then there exists a unique $\sigma \in G$\ such that $S\sigma =T$ and $%
Sw\sigma =Tw$.\ In particular, we have\ $Sw=S$ if and only if $Tw=T$.

For each $\Delta $-tiling $\mathcal{T}$, each $T\in \mathcal{T}$\ and any $%
v,w\in \Omega _{\mathcal{L}}$ such that $Tvw$\ and $Twv$\ exist, consider $%
\sigma ,\tau \in G$ such that $Tv=T\sigma $, $Twv=Tw\sigma $, $Tw=T\tau $
and $Tvw=Tv\tau $.\ Then we have $Tvw=T\sigma \tau =T\tau \sigma =Twv$\
since $\sigma \tau =\tau \sigma $.~~$\blacksquare $\bigskip

\noindent \textbf{Proof of Theorem 3.12.} By Lemma 3.12.1, the $\mathcal{L}%
_{\Delta }$-structures associated to $\Delta $-tilings are equational,
commutative, and form a strongly regular class. If a class $\mathcal{C}$\ of 
$\Delta $-tilings is defined by local rules, then, by Theorem 2.1, the same
property is true for the class $K$ consisting of the associated\ $\mathcal{L}%
_{\Delta }$-structures. The structures in $K$ are periodic if the tilings in 
$\mathcal{C}$ are periodic. Then Theorem 3.11 implies that $K$, and
therefore $\mathcal{C}$, is a finite union of isomorphism classes.~~$%
\blacksquare $\bigskip

\textbf{4. Local isomorphism, rigidness and aperiodicity.}\bigskip

We say that a tiling or a relational structure is \textit{rigid} if it has
no nontrivial automorphism. In the present section, we are interested in
characterizing tilings, and more generally relational structures, which are
locally isomorphic to rigid ones. Theorem 4.1 (respectively Corollary 4.2)
gives a characterization for relational structures (respectively tilings)
which are uniformly locally finite and satisfy the local isomorphism
property. Corollary 4.3 gives a simpler characterization concerning tilings
of the euclidean spaces $%
\mathbb{R}
^{n}$, where isomorphism is defined modulo a group of isometries.

In [7], we considered tilings of the euclidean spaces $%
\mathbb{R}
^{n}$, and isomorphism was defined up to translation. In that case, [7,
Proposition 2.4] implies that the relational structure $M$ associated to a
tiling is rigid if and only if the tiling is not invariant through any
nontrivial translation. If a connected structure $N$ is locally isomorphic
to $M$, then $N$ is associated to another tiling by [7, Corollary 5.4].
Moreover, according to [7, Proposition 5.1], the tilings associated to $M$
and $N$ are invariant through the same translations of $%
\mathbb{R}
^{n}$. It follows that $N$ is rigid if and only if $M$ is rigid.

Examples 1, 2, 3, which are given after Corollary 4.3, imply that the last
property is no longer true if we consider isomorphism modulo an arbitrary
group of isometries of an euclidean space $%
\mathbb{R}
^{n}$, or tilings of a noneuclidean space. Similarly, Example 4 illustrates
Theorem 4.1 for relational structures which are not associated to tilings.
Examples 5 and 6 are given in order to show the importance of each
hypothesis in Theorem 4.1.\bigskip

\noindent \textbf{Theorem 4.1.} Let $\mathcal{L}$ be a finite relational
language and let $M$ be a uniformly locally finite $\mathcal{L}$-structure
which satisfies the local isomorphism property. Then $M$ is locally
isomorphic to a connected rigid $\mathcal{L}$-structure if and only if, for
each $r\in 
\mathbb{N}
^{\ast }$, there exists $x\in M$ such that $(M,y)\cong (M,z)\ $implies $y=z$%
\ for $y,z\in B_{M}(x,r)$.\bigskip

\noindent \textbf{Proof.} By K\"{o}nig's lemma, the following properties are
equivalent for locally finite $\mathcal{L}$-structures:

\noindent (P) There exist $r\in 
\mathbb{N}
^{\ast }$ and, for each $x\in N$, $y\neq z$ in $B_{N}(x,r)$ such that $%
(N,y)\cong (N,z)$.

\noindent (Q) There exist $r\in 
\mathbb{N}
^{\ast }$ and, for each $x\in N$ and each $s\in 
\mathbb{N}
$, $y\neq z$ in $B_{N}(x,r)$ such that $(B_{N}(y,s),y)\cong (B_{N}(z,s),z)$.

Any $\mathcal{L}$-structure $N$ which is locally isomorphic to $M$ is
uniformly locally finite like $M$. If\ $M$ satisfies (P), and therefore
satisfies (Q), then\ $N$ satisfies (Q) since it locally isomorphic to $M$,
and therefore satisfies (P). Consequently, $N$ is not rigid.

Now, let us suppose that $M$ does not satisfy (P). First we show that there
exist a sequence $(x_{n})_{n\in 
\mathbb{N}
}$ in $M$ and two strictly increasing sequences $(r_{n})_{n\in 
\mathbb{N}
}$ and $(s_{n})_{n\in 
\mathbb{N}
}$ in $%
\mathbb{N}
$ such that, for each $n\in 
\mathbb{N}
$, $(B_{M}(x_{n},r_{n}+s_{n}),x_{n})\cong
(B_{M}(x_{n+1},r_{n}+s_{n}),x_{n+1})$ and $B_{M}(x_{n},r_{n})$ contains no
elements $y\neq z$ with $(B_{M}(y,s_{n}),y)\cong (B_{M}(z,s_{n}),z)$. We
write $r_{0}=s_{0}=0$ and we take for $x_{0}$\ any element of $M$.

For each $n\in 
\mathbb{N}
$, supposing that $x_{n},r_{n},s_{n}$\ are already defined, we define $%
x_{n+1},r_{n+1},s_{n+1}$ as follows: As $M$ satisfies the local isomorphism
property, there exists $r\in 
\mathbb{N}
$ such that each $B_{M}(x,r)$\ contains some $y$ with $%
(B_{M}(y,r_{n}+s_{n}),y)\cong (B_{M}(x_{n},r_{n}+s_{n}),x_{n})$; we take $%
r>r_{n}$. As $M$ does not satisfy (Q), there exist $x\in M$\ and $s\in 
\mathbb{N}
^{\ast }$ such that $B_{M}(x,2r)$ contains no elements $y\neq z$ with $%
(B_{M}(y,s),y)\cong (B_{M}(z,s),z)$. We take for $x_{n+1}$\ any $u\in
B_{M}(x,r)$\ such that $(B_{M}(u,r_{n}+s_{n}),u)\cong
(B_{M}(x_{n},r_{n}+s_{n}),x_{n})$. Then $B_{M}(x_{n+1},r)$ contains no
elements $y\neq z$ with $(B_{M}(y,s),y)\cong (B_{M}(z,s),z)$. We take $%
r_{n+1}=r$ and $s_{n+1}=\sup (s_{n}+1,s)$.

We consider the inductive limit $(N,x)$ of the pairs $%
(B_{M}(x_{n},r_{n}+s_{n}),x_{n})$ relative to some isomorphisms

\noindent $\theta _{n}:(B_{M}(x_{n},r_{n}+s_{n}),x_{n})\rightarrow
(B_{M}(x_{n+1},r_{n}+s_{n}),x_{n+1})$

\noindent\ $\ \ \ \ \ \ \ \ \ \ \ \ \ \ \ \ \ \ \ \ \ \ \ \ \ \ \ \ \ \ \ \
\subset (B_{M}(x_{n+1},r_{n+1}+s_{n+1}),x_{n+1})$.

\noindent As $M$ satisfies the local isomorphism property, $N$ is locally
isomorphic to $M$. For $n\in 
\mathbb{N}
$ and $y\neq z$ in $B_{N}(x,r_{n})$, we have $(B_{N}(y,s_{n}),y)\ncong
(B_{N}(z,s_{n}),z)$ since $B_{N}(y,s_{n})$ and $B_{N}(z,s_{n})$ are
contained in$\ B_{N}(x,r_{n}+s_{n})$ and $(B_{N}(x,r_{n}+s_{n}),x)$ is
isomorphic to $(B_{M}(x_{n},r_{n}+s_{n}),x_{n})$.

For each nontrivial automorphism $\theta $ of $N$, there exist $n\in 
\mathbb{N}
$ and $y\neq z$ in $B_{N}(x,r_{n})$ such that $y\theta =z$. We have $%
(B_{N}(y,s_{n}),y)\cong (B_{N}(z,s_{n}),z)$, whence a contradiction.~~$%
\blacksquare $\bigskip

From now on, we consider the metric space $(E,\delta )$, the group $G$ and
the set $\Delta $ defined in Section 1. As a consequence of Theorem 4.1, we
have:\bigskip

\noindent \textbf{Corollary 4.2.} Let $\mathcal{T}$ be a $\Delta $-tiling
which satisfies the local isomorphism property. Then $\mathcal{T}$ is
locally isomorphic to a rigid $\Delta $-tiling if and only if, for each $%
r\in 
\mathbb{N}
^{\ast }$, there exists $T\in \mathcal{T}$\ such that no $\sigma \in
G-\left\{ \mathrm{Id}\right\} $ with $\mathcal{T}\sigma =\mathcal{T}$
satisfies $T\sigma \in \mathcal{B}_{r}^{\mathcal{T}}(T)$.\bigskip

\noindent \textbf{Remark.} If $(E,\delta )$\ is weakly homogeneous, then
Corollary 4.2 and Proposition 1.2 imply that $\mathcal{T}$ is locally
isomorphic to a rigid $\Delta $-tiling if and only if, for each $\alpha \in 
\mathbb{R}
_{>0}$, there exists $x\in E$\ such that no $\sigma \in G-\left\{ \mathrm{Id}%
\right\} $ with $\mathcal{T}\sigma =\mathcal{T}$ satisfies $\delta
(x,x\sigma )\leq \alpha $.\bigskip

\noindent \textbf{Proof of Corollary 4.2.} It follows from Theorem 1.3,
Corollary 2.2 and Proposition 2.3 that $\mathcal{T}$ is locally isomorphic
to a rigid $\Delta $-tiling if and only if the associated $\mathcal{L}%
_{\Delta }$-structure is locally isomorphic to a connected rigid $\mathcal{L}%
_{\Delta }$-structure, and therefore, by Theorem 4.1, if and only if, for
each $r\in 
\mathbb{N}
^{\ast }$, there exists $T\in \mathcal{T}$ such that $(\mathcal{T},U)\cong (%
\mathcal{T},V)$ implies $U=V$ for any $U,V\in \mathcal{B}_{r}^{\mathcal{T}%
}(T)$.

Now we show that this property is true if and only if, for each $r\in 
\mathbb{N}
^{\ast }$, there exists $T\in \mathcal{T}$\ such that no $\sigma \in
G-\left\{ \mathrm{Id}\right\} $ with $\mathcal{T}\sigma =\mathcal{T}$
satisfies $T\sigma \in \mathcal{B}_{r}^{\mathcal{T}}(T)$. For each $T\in 
\mathcal{T}$, each $r\in 
\mathbb{N}
^{\ast }$ and each $\sigma \in G-\left\{ \mathrm{Id}\right\} $\ such that $%
\mathcal{T}\sigma =\mathcal{T}$, if there exist $U\neq V$ in $\mathcal{B}%
_{r}^{\mathcal{T}}(T)$ such that $U\sigma =V$, then we have $T\sigma \in 
\mathcal{B}_{3r}^{\mathcal{T}}(T)$. Conversely, if $r\geq 2q$ for the
integer $q$ of Section 1 and if $T\sigma \in \mathcal{B}_{r}^{\mathcal{T}%
}(T) $, then we obtain $U\neq V$ in $\mathcal{B}_{r}^{\mathcal{T}}(T)$ such
that $U\sigma =V$ as follows: we write $U=T$\ and $V=T\sigma $ if $T\sigma
\neq T$; otherwise, we consider any $U\in \mathcal{B}_{q}^{\mathcal{T}}(T)$\
such that $U\sigma \neq U$, and we write $V=U\sigma $.~~$\blacksquare $%
\bigskip

\noindent \textbf{Corollary 4.3.} Let $n\geq 1$\ be an integer and let $%
\mathcal{T}$ be a tiling of the euclidean space $%
\mathbb{R}
^{n}$ which satisfies the local isomorphism property. Then $\mathcal{T}$ is
locally isomorphic to a rigid tiling if and only if it is not invariant
through a nontrivial translation.\bigskip

\noindent \textbf{Proof.} By the remark after Corollary 4.2, it suffices to
show that, if $\mathcal{T}$ is not invariant through a nontrivial
translation, then, for each $\alpha \in 
\mathbb{R}
_{>0}$, there exists $x\in 
\mathbb{R}
^{n}$ such that $\left\Vert xh-x\right\Vert >\alpha $ for each nontrivial
isometry $h$ of $%
\mathbb{R}
^{n}$ which stabilizes $\mathcal{T}$. This result follows from Theorem 4.4
below since the isometries of $%
\mathbb{R}
^{n}$ which stabilize $\mathcal{T}$\ form a discrete group by Proposition
1.8.~~$\blacksquare $\bigskip

\noindent \textbf{Theorem 4.4.} Let $n\geq 1$\ be an integer and let $H$ be
a discrete group of isometries of the euclidean space $%
\mathbb{R}
^{n}$. Then $H$ contains no translation if and only if, for each $\alpha \in 
\mathbb{R}
_{>0}$, there exists $x\in 
\mathbb{R}
^{n}$ such that $\left\Vert xh-x\right\Vert >\alpha $ for each $h\in
H-\left\{ \mathrm{Id}\right\} $.\bigskip

\noindent \textbf{Proof.} The \textquotedblleft if\textquotedblright\ part
is clear. For the \textquotedblleft only if\textquotedblright\ part, we
proceed by induction on $n$. We denote by $\mathrm{E}_{n}$\ the group of all
isometries of $%
\mathbb{R}
^{n}$. For any $X,Y\subset 
\mathbb{R}
^{n}$, we write $\delta (X,Y)=\inf_{x\in X,y\in Y}\left\Vert y-x\right\Vert $%
.

For each $f\in \mathrm{E}_{n}$, we consider the linear map $\overline{f}$\
associated to $f$, the largest affine subspace $W_{f}$ of $%
\mathbb{R}
^{n}$\ with $W_{f}f=W_{f}$\ such that $f$\ acts on $W_{f}$\ as a
translation,\ and the element $w_{f}\in 
\mathbb{R}
^{n}$ such that $xf=x+w_{f}$\ for each $x\in W_{f}$.\ We have\ $W_{f}=\{x\in 
\mathbb{R}
^{n}\mid \left\Vert xf-x\right\Vert $ minimal$\}$.

If $V_{f}$ is the vector subspace of $%
\mathbb{R}
^{n}$\ orthogonal to $W_{f}$ and maximal for that property, then the
restriction of $\overline{f}$ to $V_{f}$ is an orthogonal transformation
without nontrivial fixed point. Consequently, for each $\alpha \in 
\mathbb{R}
_{>0}$, there exists $\beta \in 
\mathbb{R}
_{>0}$ such that, for each $x\in 
\mathbb{R}
^{n}$, $\delta (x,W_{f})>\beta $ implies $\left\Vert xf-x\right\Vert >\alpha 
$.

It follows from a theorem of Bieberbach (see [9, Theorem 1, p. 15]) that $H$
has a normal subgroup $N$ with $N$ abelian and $H/N$ finite.\ As $N$ is
finitely generated, the same properties remain true if we replace $N$ by $%
N^{r}=\{h^{r}\mid h\in N\}$\ for an integer $r\geq 2$. Consequently, we can
suppose for the remainder of the proof that $N$ is torsion-free. Then we
have $w_{f}\neq 0$ for each $f\in N-\left\{ \mathrm{Id}\right\} $\ since\ $N$
is discrete like $H$.

We observe that\ $W_{f}g=W_{f}$ for any $f,g\in \mathrm{E}_{n}$ which
commute, and in particular for $f,g\in N$:\ We have $%
(W_{f}g)f=(W_{f}f)g=W_{f}g$\ and

\noindent $(zg)f-(yg)f=(zf)g-(yf)g=(zf-yf)\overline{g}=(z-y)\overline{g}%
=zg-yg$

\noindent for any $y,z\in W_{f}$. It follows that $f$ stabilizes $W_{f}g$
and acts on $W_{f}g$\ as a translation, which implies $W_{f}g=W_{f}$.

We consider $W=\cap _{f\in N}W_{f}$.\ In order to prove that $W$ is
nonempty, it suffices to show that, for each $r\in 
\mathbb{N}
^{\ast }$ and any $f_{1},...,f_{r+1}\in N$, if $W_{r}=W_{f_{1}}\cap ...\cap
W_{f_{r}}$\ is nonempty, then $W_{r+1}=W_{f_{1}}\cap ...\cap W_{f_{r+1}}$\
is nonempty. But we have $W_{r}f_{r+1}=W_{r}$\ since $%
W_{f_{i}}f_{r+1}=W_{f_{i}}$\ for $1\leq i\leq r$; it follows that $W_{r+1}$
is the largest affine subspace $V$ of $W_{r}$\ with $Vf_{r+1}=V$\ such that $%
f_{r+1}$\ acts on $V$\ as a translation.

Now we show that $Wh=W$\ for each $h\in H$. For each $f\in N$, we have $%
W(hfh^{-1})=W$\ since $hfh^{-1}\in N$, and therefore $Whf=W(hfh^{-1})h=Wh$.
Moreover, for each $f\in N$ and any $y,z\in W$, we have

\noindent $zhf-yhf=z(hfh^{-1})h-y(hfh^{-1})h=[z(hfh^{-1})-y(hfh^{-1})]%
\overline{h}$

$=(z-y)\overline{h}=zh-yh$

\noindent since $hfh^{-1}\in N$. It follows that each $f\in N$ stabilizes $%
Wh $ and acts on $Wh$\ as a translation, which implies $Wh=W$.

Now we fix $\alpha \in 
\mathbb{R}
_{>0}$\ and we prove that there exists $x\in 
\mathbb{R}
^{n}$ such that $\left\Vert xh-x\right\Vert >\alpha $ for each $h\in
H-\left\{ \mathrm{Id}\right\} $. We consider the set $\Omega $ of all affine
subspaces\ of $%
\mathbb{R}
^{n}$\ which are orthogonal to $W$ and maximal for that property.\ For each $%
U\in \Omega $\ and each $h\in H$, we have $Uh\in \Omega $\ since $Wh=W$.

First we show that $\{g\in H\mid \delta (U,Ug)\leq \alpha \}$\ is finite for
each $U\in \Omega $. As each $h\in N$ is acting on $W$\ as a translation of
vector $w_{h}$, we have $\delta (V,Vh)=\left\Vert w_{h}\right\Vert $\ for $%
V\in \Omega $\ and $h\in N$, and $w_{gh}=w_{g}+w_{h}$\ for $g,h\in N$. We
write $\gamma =\inf_{h\in N-\left\{ \mathrm{Id}\right\} }\left\Vert
w_{h}\right\Vert $. We have $\gamma \neq 0$\ since $N$\ is discrete and
torsion-free. We consider $r\in 
\mathbb{N}
^{\ast }$ such that $r\gamma >2\alpha $. For each $U\in \Omega $, each $g\in
H$\ such that\ $\delta (U,Ug)\leq \alpha $ and each $h\in N-\left\{ \mathrm{%
Id}\right\} $, the inequality\ $\delta (Ug,Ugh^{r})=\left\Vert
w_{h^{r}}\right\Vert =r\left\Vert w_{h}\right\Vert \geq r\gamma >2\alpha $\
implies $\delta (U,Ugh^{r})>\alpha $. As $N^{r}$\ has finite index in $N$\
and therefore in $H$, it follows that $\{g\in H\mid \delta (U,Ug)\leq \alpha
\}$\ is finite.

Now we consider $K=\{h\in H\mid xh-x\in W$\ for each $x\in 
\mathbb{R}
^{n}\}$\ and, for each $h\in K$, the\ restriction $h_{W}$\ of $h$ to $W$.
For each $x\in 
\mathbb{R}
^{n}$, we denote by $x_{W}$ the projection of $x$ on $W$. Then $%
K_{W}=\{h_{W}\mid h\in K\}$\ is, like $K$, a discrete group of isometries\
without nontrivial translation, since $xh-x=x_{W}h-x_{W}$ for each $h\in K$
and each $x\in 
\mathbb{R}
^{n}$. The induction hypothesis applied to $W$\ and $K_{W}$\ implies that
there exists $x\in W$\ such that $\left\Vert xh_{W}-x\right\Vert >\alpha $
for each $h\in K-\left\{ \mathrm{Id}\right\} $.

We consider the unique $U\in \Omega $ such that $x\in U$. We have $\delta
(U,Uh)=\left\Vert xh_{W}-x\right\Vert >\alpha $ for each $h\in K-\left\{ 
\mathrm{Id}\right\} $. For each $h\in H-K$, we have $U\cap
W_{h}\varsubsetneq U$ and $S_{h}=\{y\in U\mid \left\Vert yh-y\right\Vert
\leq \alpha \}$ is contained in $A+(U\cap W_{h})$ for a bounded subset $A$
of $U$ since there exists $\beta \in 
\mathbb{R}
_{>0}$ such that, for each $x\in 
\mathbb{R}
^{n}$, $\delta (x,W_{h})>\beta $ implies $\left\Vert xh-x\right\Vert >\alpha 
$. As $S_{h}$\ is empty for each $h\in H$ such that $\delta (U,Uh)>\alpha $,
there exist finitely many nonempty subsets $S_{h}$, and their union cannot
cover $U$.~~$\blacksquare $\bigskip

Now, we illustrate Corollary 4.2 and Corollary 4.3 with three examples
related to aperiodicity. Several different definitions have been given for
that notion (see [4, p. 4]).

We consider the system $\Delta $ defined in Section 1 and the set $\mathcal{C%
}$\ of all $\Delta $-tilings\ which satisfy\ a set $\Omega $ of \emph{local
rules}, each of them saying which configurations of some given size can
appear in a tiling belonging to $\mathcal{C}$.\ According to [12, p. 208],
we say that $\Omega $ is \emph{strong} if the $\Delta $-tilings in $\mathcal{%
C}$ satisfy the local isomorphism property and if they are not invariant
through any nontrivial translation. We use the classical definition of
translation for the euclidean spaces $%
\mathbb{R}
^{n}$, and the definition given in Section 3 for the general case.

Here we do not suppose $\Omega $\ finite. One reason is that some natural
sets of tilings are defined by strong infinite sets of local rules (for
instance we showed this property in [8] for the set of all complete folding
sequences, and for the set of all coverings of the plane by complete folding
curves which satisfy the local isomorphism property). Another reason is
that, for each $\Delta $-tiling $\mathcal{T}$ which satisfies the local
isomorphism property, the set of all $\Delta $-tilings\ which are locally
isomorphic to $\mathcal{T}$ is defined by a set of local rules which can be
finite as in Examples 2 and 3, or infinite as in Example 1.

By Corollary 4.3, any tiling of an euclidean space of finite dimension which
satisfies a strong set of local rules is locally isomorphic to a rigid
tiling if and only if it is not invariant through a nontrivial translation.
We do not know presently if this result can be generalized with the notion
of translation that we consider.

In Examples 1 and 3 below, the group $G$ consists of all isometries of $E$;
in Example 2 we only consider positive isometries of $%
\mathbb{R}
^{3}$. In each example, all the structures are uniformly locally finite and
satisfy the local isomorphism property. On the other hand, they do not
satisfy (P). Some of them are rigid and others have nontrivial
automorphisms, but all of them are locally isomorphic.\bigskip

\noindent \textbf{Example 1.} For each $r\in 
\mathbb{R}
-%
\mathbb{Q}
$ and each $s\in 
\mathbb{R}
$, we consider the line $L(r,s)$\ of equation $y=rx+s$. We write

\noindent $\Omega (r,s)=\{(a,b)\in 
\mathbb{Z}
\times 
\mathbb{Z}
\mid L(r,s)\cap ([a-1/2,a+1/2[\times \lbrack b-1/2,b+1/2[)\neq \varnothing
\} $.

\noindent For each $a\in 
\mathbb{Z}
$, $\left\{ a\right\} \times 
\mathbb{Z}
$ contains $n+1$ or $n+2$\ points of $\Omega (r,s)$, where $n$ is the
integral part of $\left\vert r\right\vert $. We colour the point $a$ in
white if $\left\{ a\right\} \times 
\mathbb{Z}
$ contains $n+1$\ points of $\Omega (r,s)$ and in black otherwise. We
consider the tiling $\mathcal{T}(r,s)$ of the euclidean space $%
\mathbb{R}
$ which consists of the segments $[a,a+1]$\ for $a\in 
\mathbb{Z}
$ with their endpoints coloured in white or black as above.

For each $r\in 
\mathbb{R}
-%
\mathbb{Q}
$, the tilings $\mathcal{T}(r,s)$ are locally isomorphic and each of them
satisfies the local isomorphism property. They do not satisfy (P) since they
are not invariant through any nontrivial translation. They are invariant
through a unique symmetry if there exist $a,b\in 
\mathbb{Z}
$\ such that $(a,b)$\ or $(a+1/2,b)$\ or $(a,b+1/2)$\ belongs to $L(r,s)$,
and rigid otherwise. The three possibilities above are respectively realized
for $s\in 
\mathbb{Z}
+r%
\mathbb{Z}
$, $s\in r/2+%
\mathbb{Z}
+r%
\mathbb{Z}
$,\ $s\in 1/2+%
\mathbb{Z}
+r%
\mathbb{Z}
$.\bigskip

\noindent \textbf{Example 2.} Let $(E,\delta )$\ be the euclidean space $%
\mathbb{R}
^{3}$ and let $G$ consist of the positive isometries. Let \emph{T} be\ a
tiling of $%
\mathbb{R}
^{3}$ which satisfies the local isomorphism property. Suppose that the group
of isometries which leave $\mathcal{T}$\ globally invariant is generated by
a "screwing motion" $\sigma $, which is the composition of a translation
with a rotation about an axis parallel to the translation. If the angle of
the rotation belongs to $\pi 
\mathbb{Q}
$,\ then some nontrivial power of $\sigma $\ is a translation, and $\mathcal{%
T}$\ satisfies (P). Otherwise, $\mathcal{T}$ does not satisfy (P), and
Theorem 4.1 implies that $\mathcal{T}$ is locally isomorphic to a rigid
tiling. According to [12, Section 7.2, pp. 208-213], examples of that
situation have been given by Danzer for tilings obtained from $1$\ prototile
(the examples with $n$ odd must be considered).\bigskip

\noindent \textbf{Example 3.} In 1979, R. Penrose gave his famous example
(see [3]) of two polygonal prototiles, the \textquotedblleft
arrow\textquotedblright\ and the \textquotedblleft kite\textquotedblright ,
which define an aperiodic class of tilings of the euclidean space $%
\mathbb{R}
^{2}$. There exist $2^{\omega }$\ Penrose tilings. All of them are locally
isomorphic and each of them satisfies the local isomorphism property. The
Robinson tilings (see [10]) have the same properties, but they are
constructed from a larger set of prototiles. Penrose asked if there exist
classes of tilings of $%
\mathbb{R}
^{2}$ defined from a single prototile which have these properties. The
question is apparently still open for tilings with non-overlapping tiles
(see [6]).

In the hyperbolic plane, it is not difficult to construct such an example.
Here, we use the representation of the hyperbolic plane by the Poincar\'{e}
half-plane $%
\mathbb{R}
\times 
\mathbb{R}
_{>0}$.

Figure 1 illustrates the construction of the tilings in our example, which
is a particular case of those given in [6]. We denote by $\Omega $\ the set
of all tilings constructed in that way. For any $\mathcal{S},\mathcal{T}\in
\Omega $, each $S\in \mathcal{S}$\ and each $T\in \mathcal{T}$, there exists
a unique $\sigma \in G$\ such that $S\sigma =T$, because of the arrows on
the edges of the tiles.

For each $\mathcal{T}\in \Omega $\ and each $T\in \mathcal{T}$, we consider $%
(U_{n}(T))_{n\in 
\mathbb{N}
}$ where $U_{0}(T)=T$ and, for each $n\in 
\mathbb{N}
$, $U_{n+1}(T)$ is the tile just above $U_{n}(T)$. For each $n\in 
\mathbb{N}
^{\ast }$, we write $a_{n}(T)=0$ if $U_{n}(T)$ is at the left of $U_{n-1}(T)$%
, and $a_{n}(T)=1$ otherwise.

For any $\mathcal{S},\mathcal{T}\in \Omega $, each $S\in \mathcal{S}$\ and
each $T\in \mathcal{T}$, we have $(\mathcal{S},S)\cong (\mathcal{T},T)$ if
and only if $(a_{n}(S))_{n\in 
\mathbb{N}
^{\ast }}=(a_{n}(T))_{n\in 
\mathbb{N}
^{\ast }}$. For each $r\in 
\mathbb{N}
^{\ast }$, we have $(\mathcal{B}_{r}^{\mathcal{S}}(S),S)\cong (\mathcal{B}%
_{r}^{\mathcal{T}}(T),T)$ if and only if $%
(a_{1}(S),...,a_{r}(S))=(a_{1}(T),...,a_{r}(T))$.

For any $\mathcal{S},\mathcal{T}\in \Omega $, each $S\in \mathcal{S}$, each $%
T\in \mathcal{T}$ and each $r\in 
\mathbb{N}
^{\ast }$, there exists $S^{\prime }\in \mathcal{T}$\ such that $%
T=U_{r}(S^{\prime })$ and $(a_{1}(S^{\prime }),...,a_{r}(S^{\prime
}))=(a_{1}(S),...,a_{r}(S))$, which implies $S^{\prime }\in \mathcal{B}_{r}^{%
\mathcal{T}}(T)$ and $(\mathcal{B}_{r}^{\mathcal{S}}(S),S)\cong (\mathcal{B}%
_{r}^{\mathcal{T}}(S^{\prime }),S^{\prime })$. Consequently, any tiling in $%
\Omega $ satisfies the local isomorphism property, and any two such tilings
are locally isomorphic.

For each $\mathcal{T}\in \Omega $ and any $S,T\in \mathcal{T}$, there exist $%
i,j\in 
\mathbb{N}
$\ such that $U_{i}(S)=U_{j}(T)$; for $n\geq i+1$, we have $%
a_{n}(S)=a_{n+k}(T)$ where $k=j-i$. Consequently, each $\mathcal{T}\in
\Omega $ only realizes countably many sequences $(a_{n})_{n\in 
\mathbb{N}
^{\ast }}\in \{0,1\}^{%
\mathbb{N}
^{\ast }}$. As each such sequence is realized by a tiling $\mathcal{T}\in
\Omega $, it follows that $\Omega $ is the union of $2^{\omega }$\
isomorphism classes. This property is a particular case of Corollary 2.6
above.

Now we show that, for each $\mathcal{T}\in \Omega $ which is not rigid\ and
each $T\in \mathcal{T}$, there exists $k\in 
\mathbb{N}
^{\ast }$\ such that $a_{n}(T)=a_{n+k}(T)$ for $n$ large enough: We consider 
$\sigma \in G-\{\mathrm{Id}\}$ such that $\mathcal{T}\sigma =\mathcal{T}$
and $i,j\in 
\mathbb{N}
$\ such that $U_{i}(T)=U_{j}(T\sigma )$. We have $i\neq j$ because $%
U_{i}(T)=U_{i}(T\sigma )=U_{i}(T)\sigma $ would imply $\sigma =\mathrm{Id}$.
For $k=j-i$, we have $U_{n}(T)=U_{n+k}(T\sigma )$ for $n\geq i$, and
therefore $a_{n}(T)=a_{n+k}(T\sigma )=a_{n+k}(T)$ for $n\geq i+1$.

Conversely, for each $\mathcal{T}\in \Omega $\ and each $T\in \mathcal{T}$,
if $I=\{k\in 
\mathbb{Z}
\mid a_{n}(T)=a_{n+k}(T)$ for $n$ large enough$\}$\ contains a nonzero
element, then $I$\ is the ideal of $%
\mathbb{Z}
$ generated by the smallest $h\in 
\mathbb{N}
^{\ast }$ which belongs to $I$. For each $r\in 
\mathbb{N}
$\ such that $a_{n}(T)=a_{n+h}(T)$ for $n>r$, the isometry which sends $%
U_{r}(T)$ to $U_{r+h}(T)$\ generates $\{\sigma \in G\mid \mathcal{T}\sigma =%
\mathcal{T}\}$.\bigskip

\noindent \textbf{Remark.} In Example 3, similar to the case of Penrose
tilings or Robinson tilings, the class of $\Delta $-tilings that we consider
is defined by a local rule which\ describes the possible configurations of
the immediate neighbours of a tile. In the case of Penrose tilings or
Robinson tilings (see [7, p. 125]), it follows that there exists a local
rule expressed by one sentence which characterizes among the connected $%
\mathcal{L}_{\Delta }$-structures those which are associated to $\Delta $%
-tilings, because no such tiling is invariant through an infinite group of
isometries. On the other hand, in Example 3, no such rule exists since any
local rule satisfied by the $\mathcal{L}_{\Delta }$-structures associated to 
$\Delta $-tilings is also satisfied by some of their quotients.\bigskip

The following example generalizes the argument of Example 3, even though the
relational structures that we consider are not represented by
tilings:\bigskip

\noindent \textbf{Example 4.} We write $\mathcal{L}=\{P_{1},...,P_{k}\}$
where $P_{1},...,P_{k}$\ are unary functional symbols, and we consider the
nonempty $\mathcal{L}$-structures $M$ which satisfy the following properties:

\noindent 1) For each $x\in M$,\ there exists one and only one pair $(i,y)$
with $1\leq i\leq k$\ and $y\in M$ such that $yP_{i}=x$;

\noindent 2) $xP_{i_{1}}...P_{i_{r}}=x$ implies $r=0$ for $r\in 
\mathbb{N}
$, $1\leq i_{1},...,i_{r}\leq k$\ and\ $x\in M$.

\noindent Each connected such structure induces a directed tree where the
pairs $(x,y)$ of consecutive vertices are characterized by the existence of
a unique $i\in \{1,...,k\}$\ such that $xP_{i}=y$; each vertex is the origin
of $k$\ edges.

In order to apply the results of the present paper, it is convenient to
consider $P_{1},...,P_{k}$\ as binary relations. Similarly to Example 3, the
nonempty $\mathcal{L}$-structures which satisfy 1) and 2) are locally
isomorphic, and each of them satisfies the local isomorphism property. In
fact, for any such structures $M,N$, each $x\in M$, each $y\in N$\ and each $%
r\in 
\mathbb{N}
^{\ast }$, we have $(B_{M}(x,r),x)\cong (B_{M}(z,r),z)$\ for $%
z=yP_{i_{r}}...P_{i_{1}}$\ where $i_{1},...,i_{r}$ are the elements of $%
\{1,...,k\}$\ such that $xP_{i_{1}}^{-1}...P_{i_{r}}^{-1}$\ exists.

Now, for each nonempty connected $\mathcal{L}$-structure $M$ which satisfies
1), 2) and each $x\in M$, we consider the sequence $(i_{r}(x))_{r\in 
\mathbb{N}
}\in \{1,...,k\}^{%
\mathbb{N}
}$\ such that $xP_{i_{1}(x)}^{-1}...P_{i_{r}(x)}^{-1}$\ exists for each $%
r\in 
\mathbb{N}
$. Similarly to Example 3, $M$ has a nontrivial automorphism if and only if
there exists an integer $k$ such that $i_{r}(x)=i_{r+k}(x)$\ for $r$ large
enough. In that case, the group of automorphisms of $M$\ is infinite cyclic.

By Theorem 4.1, it follows that the nonempty $\mathcal{L}$-structures which
satisfy 1) and 2) are rigid.\bigskip

The last two examples are not related to tilings. They are given in order to
illustrate the importance of each hypothesis in Theorem 4.1.\bigskip

\noindent \textbf{Example 5.} The\emph{\ Cayley graph} of a group $G$
relative to a generating family $(x_{i})_{i\in I}$ is the relational
structure $M$\ defined on $G$ as follows: for $i\in I$ and $y,z\in G$, we
write $R_{i}(y,z)$ if and only if $z=yx_{i}$. The structure $M$ is uniformly
locally finite if $I$ is finite. The automorphisms of $M$ are the maps $%
y\rightarrow gy$ for $g\in G$. For any $y,z\in M$, we have $(M,y)\cong (M,z)$%
\ since there exists $g\in G$ such that $gy=z$. In particular, $M$ satisfies
the local isomorphism property and $M$ is not locally isomorphic to a rigid
structure. If $G$ is freely generated by the elements $x_{i}$, then $M$ is
not invariant through any nontrivial translation since, for each $g\in G$
and each $r\in 
\mathbb{N}
$, there exists $x\in M$ such that $gx\notin B_{M}(x,r)$.\bigskip

\noindent \textbf{Example 6.} Here, the language $\mathcal{L}$ consists of
one binary relational symbol. The $\mathcal{L}$-structure $M$ shown by
Figure 2 is uniformly locally finite, but it does not satisfy the local
isomorphism property. The only automorphisms of $M$ are the maps $%
x_{i,j}\rightarrow x_{i+k,j}$ for $k\in 
\mathbb{Z}
$. Consequently, $M$ does not satisfy the characterization of Theorem 4.1.
Anyway, each connected $\mathcal{L}$-structure $N$ locally isomorphic to $M$
is isomorphic to $M$, and therefore not rigid.\bigskip

\bigskip

\begin{center}
\textbf{References}
\end{center}

\bigskip

\noindent \lbrack 1] A. Ballier, B. Durand and E. Jeandel, Structural
aspects of tilings, in STACS 2008 (25th Annual Symposium on Theoretical
Aspects of Computer Science, Bordeaux, France, 2008), Dagstuhl Seminar
Proceedings, Internationales Begegnungs und Forschungszentrum fuer
Informatik, Schloss Dagstuhl, Germany, pp. 61-72.

\noindent \lbrack 2] Y. Benyamini and J. Lindenstrauss, Geometric nonlinear
functional analysis I, AMS Colloquium Publications 48, American Mathematical
Society, Providence, USA, 2000.

\noindent \lbrack 3] M. Gardner, Extraordinary nonperiodic tiling that
enriches the theory of tiles, Scientific American, Jan. 1977, pp. 110-121.

\noindent \lbrack 4] C. Goodman-Strauss, Open questions in tiling, notes
available at:

\noindent comp.uark.edu/\symbol{126}cgstraus/papers.

\noindent \lbrack 5] W. Hodges, Model Theory, Encyclopedia of Mathematics
and its Applications, Cambridge University Press, Cambridge, 1995.

\noindent \lbrack 6] G.A. Margulis and S. Mozes, Aperiodic tilings of the
hyperbolic plane by convex polygons, Israel J. Math. 107 (1998), 319-325.

\noindent \lbrack 7] F. Oger, Algebraic and model-theoretic properties of
tilings, Theoret. Comput. Sci. 319 (2004), 103-126.

\noindent \lbrack 8] F. Oger, Paperfolding sequences, paperfolding curves
and local isomorphism, to appear.

\noindent \lbrack 9] R.K. Oliver, On Bieberbach's analysis of discrete
euclidean groups, Proc. Amer. Math. Soc. 80 (1980), 15-21.

\noindent \lbrack 10] B. Poizat, Une th\'{e}orie finiment axiomatisable et
superstable, Groupe d'Etude de Th\'{e}ories Stables 3, Universit\'{e} Pierre
et Marie Curie, Paris, 1983.

\noindent \lbrack 11] C. Radin and M. Wolff, Space tiling and local
isomorphism, Geometriae Dedicata 42 (1992), 355-360.

\noindent \lbrack 12] M. Senechal, Quasicrystals and Geometry, Cambridge
University Press, Cambridge, GB, 1996.\bigskip

\vfill

\end{document}